\documentclass[graybox]{svmult}

\usepackage{type1cm}        % activate if the above 3 fonts are
\usepackage{makeidx}         % allows index generation
\usepackage{graphicx}        % standard LaTeX graphics tool
\usepackage{multicol}        % used for the two-column index
\usepackage[bottom]{footmisc}% places footnotes at page bottom

\usepackage{newtxtext}       % 
\usepackage{newtxmath}       % selects Times Roman as basic font

\usepackage{todonotes}

\makeindex             % used for the subject index
\usepackage{xcolor}
\usepackage{multirow}

\usepackage{float}
\usepackage{graphicx}
\usepackage{subfig}

\newcommand{\R}{\mathbb R}
\newcommand{\bs}[1]{{#1}} % use normal
\newcommand{\bu}{\bs{u}}
\newcommand{\bv}{\bs{v}}
\newcommand{\eps}{\varepsilon}
\newcommand{\Div}{\operatorname{div}}
\begin{document}

\title*{Adaptive and Pressure-Robust Discretization of Incompressible Pressure-Driven Phase-Field Fracture}
\titlerunning{Incompressible Pressure-Driven Phase-Field Fracture}
\author{Seshadri Basava and Katrin Mang and Mirjam Walloth and Thomas
  Wick and Winnifried Wollner }
\authorrunning{Basava, Mang, Walloth, Wick, Wollner}
\institute{Seshadri Basava \at Technische Universit\"at Darmstadt \email{basava@mathematik.tu-darmstadt.de}
\and Katrin Mang \at Leibniz Universit\"at Hannover \email{mang@ifam.uni-hannover.de}
\and Mirjam Walloth \at Technische Universit\"at Darmstadt \email{walloth@mathematik.tu-darmstadt.de}
\and Thomas Wick \at Leibniz Universit\"at Hannover \email{thomas.wick@ifam.uni-hannover.de}
\and Winnifried Wollner \at Technische Universit\"at Darmstadt \email{wollner@mathematik.tu-darmstadt.de}
}
%
% Use the package "url.sty" to avoid
% problems with special characters
% used in your e-mail or web address
%
\maketitle

%%%%%%%%%%%%%%%%%%%%%%%%%%%%%%%%%%%%%%%%%%%%%%%%%%%%%%
\abstract{In this work, we consider pressurized phase-field fracture 
problems in nearly and fully incompressible materials. To this end,
a mixed form for the solid equations is proposed. To enhance 
the accuracy of the spatial discretization, a residual-type error 
estimator is developed. Our algorithmic advancements are substantiated 
with several numerical tests that are inspired from benchmark configurations.
Therein, a primal-based  formulation is compared to our newly 
developed mixed phase-field fracture method for Poisson ratios 
approaching $\nu \to 0.5$. Finally, for $\nu = 0.5$, we compare the
numerical results of the mixed formulation with a pressure robust modification.
}

%%%%%%%%%%%%%%%%%%%%%%%%%%%%%%%%%%%%%%%%%%%%%%%%%%%%%%
\section{Introduction}
\label{sec:1}
This work is devoted to pressurized fractures in nearly and fully incompressible solids
using an adaptive finite element discretization.
Pressurized fracture problems modeled with a phase-field method
is currently a topic being investigated by many groups; see for instance 
\cite{WheWiWo14,BourChuYo19,Miehe2015186,HEIDER2018116,LeeWheWi16}, to name a few. 
We further extended our pressurized phase-field fracture approach
to non-isothermal configurations~\cite{NoiiWi19}.
A recent overview on pressurized and fluid-filled fractures is provided in~\cite{WheWiLee20}.
However, all these contributions deal with compressible solids in which 
Poisson's ratio is significantly less than $0.5$, i.e., the incompressible limit.

Incompressible solids are however an important field in solids mechanics~\cite{holzapfel2002nonlinear,taylor2011isogeometric,holzapfel1996large,schroder2005variational,kubo2017velocity}. 
In~\cite{MaWiWo20} a model 
and robust discretization using a phase-field method for fractures in solids 
mechanics was proposed. A well-known challenge in phase-field methods is the relationship between 
the model regularization $\eps>0$ and the spatial mesh size $h$. 
To obtain 
accurate discretizations for small $\eps$ around the fracture 
and specifically at the fracture tip adaptive 
mesh refinement is a useful tool. First studies date back to~\cite{BuOrSue10,BuOrSue13}
investigating residual-type error estimators. A predictor-corrector 
mesh refinement algorithm with a focus on crack-oriented refinement 
was developed in~\cite{HeWheWi15} and extended to three spatial dimensions in~\cite{LeeWheWi16}.
In~\cite{ArFoMiPe15}, anisotropic mesh 
refinement was studied. Goal-oriented adjoint-based a posteriori error estimation 
was subject in~\cite{Wi16_dwr_pff}. 
Based on a recent approach for residual-type a posteriori estimators for contact problems~\cite{KrauseVeeserWalloth:2015,Walloth:2019}
we developed  in~\cite{Walloth:2018} a reliable and efficient estimator for a singularly-perturbed obstacle problem taking into account the robustness (in terms of $\eps$).  We tested the resulting residual-type estimator for different fracture phase-field problems enforcing the irreversibility condition in~\cite{MangWallothWickWollner:2019} and further for nearly incompressible solids in~\cite{MaWaWiWo20}.

The main objective of the current work is two-fold. We first 
develop a phase-field model using a mixed system for pressurized fractures. 
Therein the methodology from~\cite{MaWiWo20} is combined with 
pressurized fractures as proposed in~\cite{MiWheWi19,MiWheWi15b,WheWiWo14}.
Our second aim is to apply adaptive refinement based on our residual-type error estimator~\cite{MangWallothWickWollner:2019, Walloth:2018} to this mixed-system phase-field fracture approach. 
These algorithmic concepts are substantiated with the help of several numerical examples
and mesh convergence studies
comparing classical primal formulations and our newly developed mixed formulation.
Finally, we will test a pressure-robust modification of the discrete mixed
formulation, inspired by the works~\cite{LinkeMatthiesTobiska:2016,LinkeMerdonWollner:2015} for the Stokes problem.

As this book chapter summarizes our efforts within the German Priority Programme 1748 (DFG SPP 1748),
in the project `Structure Preserving Adaptive Enriched Galerkin Methods for Pressure-Driven 3D Fracture Phase-Field Models', we briefly mention the other research directions,
which were related to our own overall goal. 

In~\cite{BrWiBeNoRa20}, we considered a stabilized decoupled iteration scheme,
a so-called $L$-scheme. Therein constant stabilization parameters 
were introduced including both numerical analysis and computational verification.
An enhancement in efficiency by using dynamically chosen stabilization 
parameters during the iteration was subsequently proposed in~\cite{EngPoWi19}.
We published our open-source 
parallel computing paper with heuristic adaptive mesh
refinement~\cite{HeiWi18_pamm}. 
The open-source programming code was used in the SPP benchmark
collection~\cite{WiMueKolletal20}.
Several comparisons of different stress-splitting methods were done 
in~\cite{FaJiWi19_paper}. 

The predictor-corrector 
approach from~\cite{HeWheWi15} inspired an adaptive non-intrusive 
global-local approach in~\cite{NoAlWiWr19}, a paper, which is also 
a collaboration within the SPP 1748 with the group of Peter Wriggers.

In the work~\cite{Walloth:2018} the basis for a provably reliable and
efficient error estimator for fracture phase-field models has been
set. The resulting residual-type error estimator has been used to
steer solely the adaptive refinement and thus the resolution of the
critical region around the crack without any prior knowledge about the
problem
in~\cite{MangWallothWickWollner:2019, MaWaWiWo20}.

The outline of this paper is as follows. 
In Section,~\ref{sec_notation_equations} 
the notation and equations are introduced. 
Next, in Section~\ref{sec_disc}, both the discretization and the numerical 
solution are addressed.
In Section~\ref{sec_error}, a residual-type error estimator 
for pressurized fractures is presented.
In the final Section~\ref{sec_tests} several numerical tests are conducted.
We summarize our findings in Section~\ref{sec_conclusions}.

%%%%%%%%%%%%%%%%%%%%%%%%%%%%%%%%%%%%%%%%%%%%%%%%%%%%%%
\section{Notation and equations}
\label{sec_notation_equations}
In this section, we introduce the basic notation and  
the underlying equations.
In the following,
let $\Omega \subset \R^2$ the total domain wherein  
$\mathcal{C}\subset \R$ denotes the fracture
and $\widehat{\Omega} \subset \R^2$ is the intact domain.
The outer boundary is denoted by $\partial \Omega$.
The inner fracture boundary is denoted by $\partial\widehat{\Omega}_F  =
\mathcal{C}$.

Using a phase-field approach, the one-dimensional fracture $\mathcal{C}$
is approximated on $\Omega\in\R^2$ 
with the help of an elliptic (Ambrosio-Tortorelli)
functional~\cite{AmTo90,AmTo92}.
This yields an approximate inner fracture boundary
$\partial\Omega_F  \approx \mathcal{C}$. 
For fracture formulations posed in a variational setting,
this has been first proposed in~\cite{BourFraMar00} based on the model 
developed in~\cite{FraMar98}.
Finally, we denote the $L^2$ scalar product with $(\cdot, \cdot)$ as 
frequently used in the literature.

Variational phase-field fracture starts with an energy functional and 
the motion of the body under consideration is then determined by the 
Euler-Lagrange equations, which are obtained by differentiation with respect 
to the unknowns.
Therefore, in phase-field-based fracture propagation, the unknown solution variables 
are vector-valued displacements $\bu:\Omega\to\mathbb{R}^2$ 
and a smoothed scalar-valued indicator phase-field function $\varphi: \Omega\to [0,1]$.
Here $\varphi = 0$ denotes the crack region
and $\varphi = 1$ characterizes the unbroken 
material. The intermediate values constitute a smooth transition zone 
dependent on a regularization parameter $\eps$.
The physics of the underlying problem ask to enforce 
a crack irreversibility condition (the crack can never heal) yielding 
the inequality constraint
\[
\varphi \leq \varphi^{n-1}.
\]
Here, $\varphi^{n-1}$ denote the previous time step solution and
$\varphi$ the current solution. 

\subsection{Pressurized phase-field fracture in a displacement formulation}
In this work, we are specifically interested in pressurized fractures in which 
a given pressure acts on the fracture boundary $\partial\Omega_F$.
Using classical interface coupling conditions, namely kinematic and dynamic coupling conditions,
for the pressure and balance of contact forces, a pressure $p_g$ can be prescribed. However,
due to the smeared zone of size $\eps$ in which $0<\varphi<1$, the exact location 
of the fracture interface is not known and leaves some freedom where to put it. 
In~\cite{MiWheWi15b}[Section 2] or~\cite{MiWheWi19}[Section 3.2], we used 
the divergence theorem to transform $p_g$ from $\partial\Omega_F$ into the 
entire domain $\Omega$. This procedure avoids knowledge of the exact fracture boundary location, but is mathematically rigorous.
Mathematical analysis (\cite{MiWheWi15b},\cite{MiWheWi19}) 
and numerous computations, 
e.g., in~\cite{WheWiWo14,HeiWi18_pamm,WiMueKolletal20}, have shown 
that this approach is justified. As a consequence of the transformation,
the pressure $p_g:\Omega\to\mathbb{R}$ is added as
domain integral to the Euler-Lagrange equations. 

Let $V:=H^1_0(\Omega; \R^2)$ and $W:=H^1(\Omega)$ the usual 
Hilbert spaces and the convex set
\[
K:= K^n = \{w\in W |\, w\leq \varphi^{n-1} \leq 1 \text{ a.e. on }
\Omega\}
\] 
including the inequality constraint.

The Euler-Lagrange system for pressurized phase-field fracture reads 
\cite{MiWheWi19}:
\begin{problem}
\label{form_1}
Let $p_g\in W^{1,\infty}(\Omega)$ be given. For the loading steps $n=1,2,3,\ldots, N$: 
Find vector-valued displacements and a
scalar-valued phase-field variable $\{\bu,\varphi\} :=
\{\bu^n,\varphi^n\} \in V \times W$ such that
\begin{equation}\label{eq_u}
\begin{aligned}
  &\Bigl(g(\varphi)\;\sigma(\bu)\,, e( {\bv}
  )\Bigr)  
  +({\varphi}^{2} p_g, \Div {\bv})
  +({\varphi}^{2} \nabla p_g, {\bv})
=0 \quad \forall \bv\in V ,
\end{aligned}
\end{equation}
and
%%%%
\begin{equation} \label{eq_varphi}
\begin{aligned}
 (1-\kappa) &({\varphi} \;\sigma(\bu):e( \bu)\,, \psi {-\varphi}) \\
&+  2 ({\varphi}\;  p_g\; \Div  \bu,\psi{-\varphi})
+ 2\, ({\varphi} \nabla p_g\cdot {\bu},\psi{-\varphi}) 
\\
&+  G_c  \Bigl( -\frac{1}{\eps} (1-\varphi,\psi{-\varphi}) + \eps (\nabla
\varphi, \nabla (\psi - {\varphi}))   \Bigr)  \geq  0
\quad \forall \psi \in K.
\end{aligned}
\end{equation}
\end{problem}
Here, 
\[
g(\varphi) = \big( (1-\kappa) {\varphi}^2  +\kappa \big)
\]
is the so-called degradation function with a small regularization parameter $\kappa$, 
$G_c$ is the critical energy release rate, and
we use the well-known Hook's law for the linear stress-strain
relationship of isotropic materials:
\begin{equation}
\label{eq_sigma_stress}
\sigma ( {\bu}) := 2\mu\, e( {\bu}) + \lambda \operatorname{tr}e( {\bu})\operatorname{I},
\end{equation}
where $\mu$ and $\lambda$ denote the Lam\'e coefficients, 
$e( {\bu}) = \frac{1}{2}(\nabla {\bu} + \nabla {\bu}^T)$ 
is the linearized strain tensor 
and $\operatorname{I}$ is the identity matrix.

\subsection{Pressurized phase-field fracture in a mixed formulation}
Following~\cite{MaWiWo20}, we now derive a mixed formulation 
for pressurized fractures.
To this end, we need to split the stress tensor~\eqref{eq_sigma_stress} into the shear 
part and the volumetric part.
In nearly incompressible materials with Poisson's ratio 
going to $0.5$, for the volumetric parameter, it holds
\[
\lambda \to \infty.
\]
To cope with volumetric locking, one possibility 
is to introduce a Lagrange multiplier, e.g.,~\cite{Br07},
with $p\in P := L^2(\Omega)$ such that
\[
p:= \lambda\; \operatorname{tr}e( {\bu}).
\]
\begin{remark}
This solution variable $p$ should not be confused with the given pressure $p_g$ from 
before. 
\end{remark}
With that, we obtain for the stress tensor:
\begin{equation*}
\sigma (\bu, p) := 2\mu\, e( {\bu}) + p\operatorname{I},
\end{equation*}
as it has been analyzed in our work~\cite{MaWiWo20} without the given
pressure $p_g$.
Adding this fracture pressure $p_g$, we obtain the following reformulation:
\begin{problem}
\label{form_1_mixed}
Let $p_g\in W^{1,\infty}(\Omega)$ be given. 
For the loading steps $n=1,2,3,\ldots, N$: Find vector-valued displacements,
a scalar-valued pressure, 
and a
scalar-valued phase-field variable 
$\{\bu,p,\varphi\} := \{\bu^n,p^n,\varphi^n\} \in V \times P \times W$ such that
\begin{equation}\label{eq_u_mixed}
\begin{aligned}
\Bigl(g(\varphi)\;\sigma(\bu,p)\,,e( {\bv})\Bigr) 
+({\varphi}^{2} p_g, \Div  {\bv}) 
+ ({\varphi}^{2} \nabla p_g,  {\bv}) 
&=0 \quad \forall \bv\in V ,
\end{aligned}
\end{equation}
and 
\begin{equation}\label{eq_p_mixed}
\begin{aligned}
(\operatorname{tr}e( {\bu})\,, q) - \frac{1}{\lambda} (p,q)
=0 \quad \forall q\in P,
\end{aligned}
\end{equation}
and
%%%%
\begin{equation} \label{eq_varphi_mixed}
\begin{aligned}
 (1-\kappa) &({\varphi} \;\sigma(\bu,p):e( \bu)\,, \psi {-\varphi}) \\
&+  2 ({\varphi}\;  p_g\; \Div  \bu,\psi{-\varphi})
+ 2\, ({\varphi} \nabla p_g\cdot  {\bu},\psi{-\varphi}) 
\\
&+  G_c  \Bigl( -\frac{1}{\eps} (1-\varphi,\psi{-\varphi}) + \eps (\nabla
\varphi, \nabla (\psi - {\varphi}))   \Bigr)  \geq  0
\quad \forall \psi \in K.
\end{aligned}
\end{equation}
\end{problem}

%%%%%%%%%%%%%%%%%%%%%%%%%%%%%%%%%%%%%%%%%%%%%%%%%%%%%%
\section{Discrete formulation}\label{sec_disc}

As the structure remains the same for all time steps, we consider one
time step $n$ for simplicity.
For the discretization in space, we decompose
the polygonal domain $\Omega$ by a (family of) meshes $\mathcal{M}=\mathcal{M}^n$ consisting
of shape regular rectangles $\mathfrak{e}$, 
such that all meshes share a
common coarse mesh. To allow for local
refinement, in particular of rectangular elements, we allow for one hanging node per edge at which degrees of freedom
will be eliminated to assert $H^1$-conformity of the discrete
spaces. 
To each mesh, we associate the mesh size function $h$, i.e.,
$h\lvert_{\mathfrak{e}} = h_{\mathfrak{e}} = \operatorname{diam}{\mathfrak{e}}$ for any
element $\mathfrak{e} \in \mathcal{M}$.
The set of nodes $q$ is given by $\mathcal{N}$ and we distinguish between the set
$\mathcal{N}^{\Gamma}$ of nodes at the boundary
and the set of interior nodes $\mathcal{N}^I$.
Later on, for the derivation of the estimator, we need the following definitions. 

For a point $q \in \mathcal{N}$, we define a patch $\omega_q$ as the interior of the union of all elements sharing the node $q$. We call the union of all sides in the interior of $\omega_q$, not including the boundary of $\omega_q$, skeleton and denote it by $\gamma_q^I$. For boundary nodes, we denote the intersections between $\Gamma$ and $\partial\omega_q$ by  $\gamma_q^{\Gamma}:=\Gamma \cap\partial\omega_q$.
Further, we will make use of $\omega_{\mathfrak{s}}$ which is the union of all elements sharing a side $\mathfrak{s}$.
We need the definition of the jump term $[\nabla\psi_h]:= \nabla|_{\mathfrak{e}}\psi_h\cdot\bs{n}_{\mathfrak{e}}- \nabla|_{\tilde{\mathfrak{e}}}\psi_h\cdot\bs{n}_{\mathfrak{e}}$ where $\mathfrak{e}, \tilde{\mathfrak{e}}$ are neighboring elements and $\bs{n}_{\mathfrak{e}}$ is the unit outward normal on the common side of the two elements. 
For the discretization, we consider (bi)-linear
($\mathbb{Q}_1(\mathfrak{e})$), (bi)-quadratic
($\mathbb{Q}_2(\mathfrak{e})$) and linear ($\mathbb{P}_1(\mathfrak{e})$) shape functions.
Thus, the finite element spaces are given by
\begin{equation*}
  \begin{aligned}
  W_h&:= W_h^n = \{\bv_h \in \mathcal{C}^0(\overline{\Omega})\mid\forall\mathfrak{e}\in\mathcal{M}, \;\bv_h|_{\mathfrak{e}}\in \mathbb{Q}_1(\mathfrak{e})\} \subset W,\\[3pt]
   P_h&:= P_h^n = \{p_h\in P \mid \forall\mathfrak{e}\in\mathcal{M}, \;p_h|_{\mathfrak{e}}\in \mathbb{P}_1(\mathfrak{e})\} \subset P,
  \end{aligned}
\end{equation*}
and  by
\begin{equation*}
 V_h:= V_h^n  = \{v_h \in \mathcal{C}^0(\overline{\Omega};\mathbb{R}^2)\mid\forall\mathfrak{e}\in\mathcal{M}, \;v_h|_{\mathfrak{e}}\in \mathbb{Q}_1(\mathfrak{e})^2 \mbox{ and }v_h=0\;\mbox{on }\Gamma \} \subset V
\end{equation*}
 for the discrete analog of Problem~\ref{form_1} 
and by 
\begin{equation*}
V_h:= V_h^n = \{v_h \in \mathcal{C}^0(\overline{\Omega};\mathbb{R}^2)\mid\forall\mathfrak{e}\in\mathcal{M}, \;v_h|_{\mathfrak{e}}\in \mathbb{Q}_2(\mathfrak{e})^2 \mbox{ and }v_h=0\;\mbox{on }\Gamma \} \subset V,
 \end{equation*}
for the discrete analogon of Problem~\ref{form_1_mixed}, respectively.

We define the respective nodal interpolation operators as $I_h^n$, and define the discrete feasible set for the phase-field by 
\begin{equation*}
K_h:=K_h^n = \{\psi_h\in W_h\mid \psi_h(q)\leq (I_h^n\varphi_h^{n-1})(q), \quad \forall q \in \mathcal{N}\}.
\end{equation*}
The nodal basis functions of the finite element space $W_h$ are denoted by $\phi_q$. 

Analogous to Problem~\ref{form_1}, we define the spatially discretized
time step problem:
\begin{problem}[Discrete formulation of Problem~\ref{form_1}]
\label{form_1_h}
Let $p_g\in W^{1,\infty}(\Omega)$ be given. For the loading steps $n=1,2,3,\ldots, N$: 
Find vector-valued displacements and a
scalar-valued phase-field variable $\{\bu_h,\varphi_h\} :=
\{\bu^n_h,\varphi^n_h\} \in V_h \times W_h$ such that
\begin{equation}\label{eq_u_disc}
\begin{aligned}
  &\Bigl(g(\varphi_h)\;\sigma(\bu_h)\,, e( {\bv_h}
  )\Bigr)  
  +({\varphi}_h^{2} p_g, \Div {\bv_h})
  +({\varphi}_h^{2} \nabla p_g, {\bv_h})
=0 \quad \forall \bv_h\in V_h ,
\end{aligned}
\end{equation}
and
%%%%
\begin{equation} \label{eq_varphi_disc}
\begin{aligned}
 (1-\kappa)& ({\varphi_h} \;\sigma(\bu_h):e( \bu_h)\,, \psi_h {-\varphi_h}) \\
&+  2 ({\varphi_h}\;  p_g\; \Div  \bu_h,\psi_h{-\varphi_h})
+ 2\, ({\varphi_h} \nabla p_g\cdot {\bu_h},\psi_h{-\varphi_h}) 
\\
&+  G_c  \Bigl( -\frac{1}{\eps} (1-\varphi_h,\psi_h{-\varphi_h}) + \eps (\nabla
\varphi_h, \nabla (\psi_h - {\varphi_h}))   \Bigr)  \geq  0
\quad \forall \psi \in K.
\end{aligned}
\end{equation}
\end{problem}

Analogous to Problem~\ref{form_1_mixed}, we define the spatially discretized
mixed time step problem:
\begin{problem}[Discrete formulation of Problem~\ref{form_1_mixed}]\label{form_1_mixed_h}
Let $p_g\in W^{1,\infty}(\Omega)$ be given. 
For the loading steps $n=1,2,3,\ldots, N$: Find vector-valued displacements,
a scalar-valued pressure, 
and a
scalar-valued phase-field variable 
$\{\bu_h,p_h,\varphi_h\} := \{\bu_h^n,p_h^n,\varphi_h^n\} \in V_h \times P_h \times W_h$ such that
\begin{equation}\label{eq_u_mixed_disc}
\begin{aligned}
\Bigl(g(\varphi_h)\;\sigma(\bu_h,p_h)\,, e( \bv_h)\Bigr)  
&+({\varphi_h}^{2} p_g, \Div  \bv_h) \\
&+ ({\varphi_h}^{2} \nabla p_g,  \bv_h) 
=0 \quad \forall \bv_h\in V_h ,
\end{aligned}
\end{equation}
and 
\begin{equation}\label{eq_p_mixed_disc}
\begin{aligned}
(\operatorname{tr}e(\bu_h)\,, q_h) - \frac{1}{\lambda} (p_h,q_h)
=0 \quad \forall q_h\in P_h,
\end{aligned}
\end{equation}
and
%%%%
\begin{equation} \label{eq_varphi_mixed_disc}
\begin{aligned}
 (1&-\kappa) (\varphi_h \;\sigma(\bu_h,p_h):e( \bu_h)\,, \psi_h -\varphi_h) \\
&+  2 (\varphi_h\;  p_g\; \Div  \bu_h,\psi_h-\varphi_h)
+ 2\, (\varphi_h \nabla p_g\cdot  \bu_h,\psi_h-\varphi_h) 
\\
&+  G_c  \Bigl( -\frac{1}{\eps} (1-\varphi_h,\psi_h-\varphi_h) + \eps (\nabla
\varphi_h, \nabla (\psi_h - \varphi_h))   \Bigr)  \geq  0
\quad \forall \psi_h \in K_h.
\end{aligned}
\end{equation}
\end{problem}
Finally, following the work of~\cite{LinkeMatthiesTobiska:2016,LinkeMerdonWollner:2015}, we propose a
pressure robust modification of Problem~\ref{form_1_mixed_h}.
To this end, we define the
divergence conforming space of Raviart-Thomas finite elements, see, e.g.,~\cite[Section~III.3.2]{BrezziFortin:1991},
on the unit square $(-1,1)^2$ by
\[
\operatorname{\mathbb{R}\mathbb{T}}_{1} = \mathbb{Q}_1^2 + \vec{x} \mathbb{Q}_1.
\]
As usual, for elements $\mathfrak{e} \in \mathcal{M}$, the space
\[
\operatorname{\mathbb{R}\mathbb{T}}_{1}(\mathfrak{e})
\]
%To this end, we define the
%divergence conforming space of Brezzi-Douglas-Marini finite elements, see, e.g.,~\cite[Section~III.3.2]{BrezziFortin:1991},
%on the unit square $(-1,1)^2$ by
%\[
%\operatorname{BDM}_{2} = P_2^2 + \operatorname{span}\Bigl(\nabla \times (x_1^3x_2), \nabla \times (x_1x_2^3)\Bigr).
%\]
%As usual, for elements $\mathfrak{e} \in \mathcal{M}$, the space
%\[
%\operatorname{BDM}_{2}(\mathfrak{e})
%\]
is then obtained by mapping of the shape functions utilizing a Piola transform.
With this, we can define the global space
\[
 \widehat{V}_h = \{v_h \in \mathcal{C}^0(\overline{\Omega};\mathbb{R}^2)\mid\forall\mathfrak{e}\in\mathcal{M}, \;v_h|_{\mathfrak{e}}\in \operatorname{\mathbb{R}\mathbb{T}}_1(\mathfrak{e})\}
 \]
 together with the interpolation operator $I_{\operatorname{RT}} \colon V_h \rightarrow \widehat{V}_h$.
 Now, following~\cite{LinkeMatthiesTobiska:2016,LinkeMerdonWollner:2015}, the pressure robust reformulation
 of Problem~\ref{form_1_mixed_h} is the problem
 \begin{problem}[Pressure robust formulation of Problem~\ref{form_1_mixed_h}]\label{form_1_mixed_h_rob}
Let $p_g\in W^{1,\infty}(\Omega)$ be given. 
For the loading steps $n=1,2,3,\ldots, N$: Find vector-valued displacements,
a scalar-valued pressure, 
and a
scalar-valued phase-field variable 
$\{\bu_h,p_h,\varphi_h\} := \{\bu_h^n,p_h^n,\varphi_h^n\} \in V_h \times P_h \times W_h$ such that
\begin{equation}\label{eq_u_mixed_disc_rob}
\begin{aligned}
\Bigl(g(\varphi_h)\;\sigma(\bu_h,p_h)\,, e( \bv_h)\Bigr)  
&+({\varphi_h}^{2} p_g, \Div  I_{\operatorname{RT}}\,\bv_h) \\
&+ ({\varphi_h}^{2} \nabla p_g,  I_{\operatorname{RT}}\,\bv_h) 
=0 \quad \forall \bv_h\in V_h ,
\end{aligned}
\end{equation}
as well as~\eqref{eq_p_mixed_disc} and~\eqref{eq_varphi_mixed_disc} hold.
 \end{problem}
 
\section{Residual-type a posteriori error estimator}
\label{sec_error}

We propose an estimator for the phase-field inequality~\eqref{eq_varphi_disc} or~\eqref{eq_varphi_mixed_disc}, respectively, to obtain a good resolution of the fracture growth. 

Utilizing either $\sigma^n_h := \sigma(\bu^n_h,p^n_h)$ for the mixed form or
$\sigma^n_h := \sigma(\bu^n_h)$ for the non-mixed form, we introduce the bilinear form
\begin{equation}\label{EpsBilinearFormDiscrete}
\begin{aligned}
a_{h,\epsilon}(\zeta,\psi) :=
&\,\frac{G_c}{\epsilon}(\zeta,\psi)  +
(1-\kappa)(\sigma^n_h:e(u^n_h)\;\zeta,\psi)  \\
&+ 2 ( p_g\, \Div  u^n_h\; \zeta,\psi)
+ 2\, ( \nabla p_g\cdot  u^n_h\;\zeta ,\psi)  + G_c\epsilon (\nabla \zeta,\nabla \psi ).
\end{aligned}
\end{equation}
Thus, the discretized variational inequality in a time step $n$ is given by
\begin{problem}[Discrete variational inequality]\label{AuxiliaryProblemDisc}
Let $u^n_h, p^n_h$ and $\varphi^{n-1}_h$ be given, then find $\varphi_h\in K_h$ such that
\begin{equation}\label{VI_DiscreteBilinearForm_Discrete}
a_{h,\epsilon}(\varphi_h , \psi_h -\varphi_h) \ge\frac{G_c}{\epsilon}\;(1 , \psi_h - \varphi_h)\quad \forall \psi_h \in K_h. 
\end{equation}
\end{problem}
We define the discrete constraining force density 
$\Lambda_{h}\in W_h^*$  of Problem~\ref{AuxiliaryProblemDisc} as
\begin{equation}\label{DiscreteConstrainingForce}
\left< \Lambda_{h},\psi_{h} \right>_{-1,1} := \frac{G_c}{\epsilon}\;(1,\psi_{h})-a_{h,\epsilon}(\varphi_{h}, \psi_{h})\quad\forall \psi_{h}\in W_h.
\end{equation}
The solution of Problem~\ref{AuxiliaryProblemDisc} is the discrete approximation of the auxiliary problem:
\begin{problem}\label{AuxiliaryProblemCont}
  Let $u^n_{h}, p^n_h$ and $\varphi^{n-1}_{h}$ be given, then find
  \[
    \hat{\varphi}\in K(I^n_{h}(\varphi^{n-1}_{h})):=\{\psi\in W\mid
    \psi\leq I^n_{h}(\varphi_{h}^{n-1})\}
  \]
  such that
\begin{equation}\label{VI_DiscreteBilinearForm}
a_{h,\epsilon}(\hat{\varphi} , \psi -\hat{\varphi} ) \ge \frac{G_c}{\epsilon}\;( 1, \psi - \hat{\varphi})\quad \forall \psi \in K( I_{h}^n(\varphi_{h}^{n-1})).
\end{equation}
\end{problem}
The corresponding constraining force density $\hat{\Lambda}\in W^*$  of Problem~\ref{AuxiliaryProblemCont} is
\begin{equation*}
\left<\hat{\Lambda} ,\psi\right>_{-1,1}:=\frac{G_c}{\epsilon}\;(1, \psi ) -a_{h,\epsilon}(\hat{\varphi}, \psi )\quad \forall \psi\in W.
\end{equation*}
\begin{remark}
As the bilinear form $a_{h,\epsilon}(\cdot,\cdot)$ depends on the approximation $u^n_h$ of $u^n$ and $p^n_h$ of $p^n$ and the constraints depend on the approximation $I^n_h(\varphi_h^{n-1})$  of $\varphi^{n-1}$, the solution $\hat{\varphi}$ of~\eqref{VI_DiscreteBilinearForm} 
is an approximation to the solution $\varphi^n$ of~\eqref{eq_varphi}
or~\eqref{eq_varphi_mixed}, respectively. 
\end{remark}
Assuming we knew $\hat{\Lambda} $ then 
\begin{align*}
\left<R(\varphi_h)\,,\psi \right>_{-1,1}:=\left<-\hat{\Lambda},\psi\right>_{-1,1} + \frac{G_c}{\epsilon}\;(1, \psi )- a_{h,\epsilon}(\varphi_h,\psi)
\end{align*}
defines the linear residual to the corresponding equation.
Thus, $R(\varphi_h)=0$ if and only if $\varphi_h=\hat{\varphi}$. Further, we are interested in the error in the constraining forces. 
As $\Lambda_h$ is not a functional on $W$, but a functional on $W_h$, it is not uniquely defined how $\Lambda_h$ acts on $W$. Thus, to compare the constraining force $\hat{\Lambda}\in W^{\ast}$ with a discrete counterpart, we choose a functional on $W^{\ast}$  called quasi-discrete constraining force, denoted by $\widetilde{\Lambda}_h\in W^*$. Therefore, we follow the approach used in~\cite{Fierro_Veeser_2003, MoonNochettoPetersdorffZhang:2007,KrauseVeeserWalloth:2015,Walloth:2018} and distinguish between full-contact nodes $q\in\mathcal{N}^{fC}$ and semi-contact nodes $q\in\mathcal{N}^{sC}$. Full-contact nodes are those nodes for which the solution is fixed to the obstacle $\varphi_h=I_h^n(\varphi^{n-1}_h)$ on $\omega_q$ and the sign condition $\left<\Lambda_h,\psi\right>_{-1,1,\omega_q}\ge0 $ $\forall \psi\ge 0\in H^1_0(\omega_q)$  is fulfilled. Semi-contact nodes are those nodes for which $\varphi_h(q)=I_h^n(\varphi^{n-1}_h)(q)$ holds but not the conditions of full-contact. 
Based on this classification, we define the quasi-discrete constraining force, where $\phi_q$ denotes the nodal basis of $W_h$, 
\begin{align}\label{QuasiDiscreteConstrainingForce}
&\left< \widetilde{\Lambda}_h,\psi \right>_{-1,1} := \sum_{q\in\mathcal{N}^{sC}}\left< \widetilde{\Lambda}_h^q,\psi\phi_q \right>_{-1,1} + \sum_{q\in\mathcal{N}^{fC}}\left< \widetilde{\Lambda}_h^q,\psi\phi_q \right>_{-1,1},
\end{align}
with the local contributions 
which are for full-contact nodes 
\begin{align*}
\left< \widetilde{\Lambda}_h^q,\psi\phi_q \right>_{-1,1}
 :=  \left<\Lambda_h, \psi\phi_q\right>_{-1,1}
\end{align*}
and for semi-contact nodes 
\begin{align*}
\left< \widetilde{\Lambda}_h^q,\psi\phi_q \right>_{-1,1}
:= \left<\Lambda_h ,\phi_q \right>_{-1,1}c_q(\psi)
\end{align*} 
with $c_q(\psi)=
\frac{\int_{\tilde{\omega}_q}\psi\phi_q}{\int_{\tilde{\omega}_q}\phi_q}$,
where $\tilde{\omega}_q$ is a proper subset of $\omega_q$.
Therefore, we define the so-called Galerkin functional
\begin{align*}
\left<G,\psi  \right>_{-1,1}:=& \left<R(\varphi_h)\,,\psi \right>_{-1,1} + \left<\hat{\Lambda}-\widetilde{\Lambda}_h,\psi \right>_{-1,1}\\
=& \Bigl(\frac{G_c}{\epsilon}, \psi \Bigr) - \left< \widetilde{\Lambda}_h, \psi\right>_{-1,1} -a_{h,\epsilon}(\varphi_h, \psi ).
\end{align*}
We note that in the case that $p_g=\mathrm{const}$ and
$\Div(u^n_h)= 0$, i.e., the material is
incompressible, the bilinear form $a_{h,\epsilon}(\zeta,\psi) $ defined in~\eqref{EpsBilinearFormDiscrete}  is elliptic; and the corresponding energy norm is given by 
\begin{equation}\label{EnergyNorm}
\|\cdot\|_{\epsilon}:= \left\{G_c\epsilon\|\nabla (\cdot)\|^2 + \|\left(\frac{G_c}{\epsilon}+(1-\kappa)\sigma(u^n_h):e(u^n_h)\right)^{\frac{1}{2}}(\cdot)\|^2 \right\}^{\frac{1}{2}}.
\end{equation}
We denote the corresponding dual norm by $\|\cdot\|_{\ast,\epsilon}:= \frac{\mathrm{sup}_{\psi\in W}\left<\cdot, \psi \right>_{-1,1}}{\|\psi\|_{\epsilon}}$.

For the definition of the error estimator contributions, we use the abbreviation of the interior residual
\begin{equation}
  \begin{aligned}
r(\varphi_h) := &\;\frac{G_c}{\epsilon} + G_c\epsilon\Delta \varphi_h -
\frac{G_c}{\epsilon}\varphi_h- (1-\kappa)(\sigma^n_h:e(u^n_h))\varphi_h \\
&+ 2 p_g \Div u^n_h \varphi_h +
2\nabla p_g\cdot u^n_h \varphi_h
\end{aligned}
\end{equation}
and set
\begin{equation}\label{Def_alpha_p}
\alpha_q:= \mathrm{min}_{x\in\omega_q}\{\frac{G_c}{\epsilon}+(1-\kappa) (\sigma(u^n_h):e(u^n_h))\}.
\end{equation}
Deriving an upper bound of $\|G\|_{\ast,\epsilon}$ as, e.g., in~\cite{Walloth:2018}, we end up with the error indicator~$\eta$ which is the sum of the following contributions
\begin{align}\label{Estimator_1}
\eta_{1}^2:=&\sum_{q\in\mathcal{N}\backslash\mathcal{N}^{fC}}\eta^2_{1,q}, & \eta_{1,q}:=&\mathrm{min}\{\frac{h_q}{\sqrt{G_c\epsilon}},\alpha_q^{-\frac{1}{2}}\}\|r(\varphi_h)\|_{\omega_q}\\ \label{Estimator_2}
\eta_{2}^2:=&\sum_{q\in\mathcal{N}\backslash\mathcal{N}^{fC}}\eta^2_{2,q}, &\eta_{2,q}:=& \mathrm{min}\{\frac{h_q}{\sqrt{G_c \epsilon}},\alpha_q^{-\frac{1}{2}}\}^{\frac{1}{2}}(G_c\epsilon)^{-\frac{1}{4}}\|G_c\epsilon[\nabla \varphi_h]\|_{\gamma_q^I}\\ \label{Estimator_3}
\eta_{3}^2:=&\sum_{q\in\mathcal{N}\backslash\mathcal{N}^{fC}}\eta^2_{3,q}, &\eta_{3,q}:=&\mathrm{min}\{\frac{h_q}{\sqrt{G_c\epsilon}},\alpha_q^{-\frac{1}{2}}\}^{\frac{1}{2}}(G_c\epsilon)^{-\frac{1}{4}}\|G_c\epsilon\nabla \varphi_h\|_{\gamma_q^N} \end{align}
In the case that $p_g=\mathrm{const}$ and the material is incompressible $\Div(u^n_h)= 0$, we can derive a robust upper bound of the error measure
\begin{align}\label{ErrorMeasure}
\|\hat{\varphi}-\varphi_h\|_{\epsilon} + \|\hat{\Lambda}-\widetilde{\Lambda}_h\|_{\ast,\epsilon}
\end{align}
in terms of the estimator
\begin{equation}\label{Def_Estimator}
\eta := \sum_{k=1}^4\eta_{k}
\end{equation}
which consists of the estimator contributions~\eqref{Estimator_1},~\eqref{Estimator_2},~\eqref{Estimator_3} and 
\begin{align*}
\eta_{4}^2:=&\sum_{q\in\mathcal{N}^{sC}}\eta^2_{4,q}, &\eta_{4,q}^2:=&s_q \int_{\widetilde{\omega}_q}(I_h^n(\varphi^{n-1}_h)-\varphi_h)\phi_q.
\end{align*}

\begin{theorem}[Reliability]
Assuming that $p_g=\mathrm{const}$ and $\Div(u^n_h)= 0$,
the error estimator $\eta$ provides a robust upper bound of the error measure, i.e.  
\begin{equation*}
\|\hat{\varphi}-\varphi_h\|_{\epsilon} + \|\hat{\Lambda}-\widetilde{\Lambda}_h\|_{\ast,\epsilon}\leq C \eta
\end{equation*}
otherwise the estimator constitutes an upper bound of the dual norm of the Galerkin functional 
$$\|G\|_{\ast,\epsilon}\leq C \eta,$$
where $C$ does not depend on $\epsilon$.
\end{theorem}

If $p_g=\mathrm{const}$ and $\Div(u^n_h)= 0$, the local estimator contributions constitute local lower bounds with respect to the local error measure~\eqref{ErrorMeasure}. 
The proof to show reliability as well as efficiency follows the ideas of~\cite{Walloth:2018}.

%%%%%%%%%%%%%%%%%%%%%%%%%%%%%%%%%%%%%%%%%%%%%%%%%%%%%%
\section{Numerical tests}
\label{sec_tests}
In this section, we investigate some examples all motivated by the theoretical calculations of Sneddon~\cite{sneddon1946distribution} and Sneddon and Lowengrub~\cite{SneddLow69} considering a pressure-driven cavity.

\begin{figure}
\begin{minipage}[b]{0.4\textwidth}
\begin{tikzpicture}
\draw[gray] (0,0)  -- (0,4) -- (4,4)  -- (4,0) -- cycle;
\node[gray] at (-0.35,0.5) {$\partial\Omega$};
\node at (3.5,3.5) {$\Omega$};
\draw [fill,opacity=0.1, blue] plot [smooth] coordinates { (1.7,2) (1.8,1.9) (2.2,1.9) (2.3,2) (2.2,2.1) (1.8,2.1) (1.7,2)};
\draw [thick, blue] plot [smooth] coordinates { (1.7,2) (1.8,1.9) (2.2,1.9) (2.3,2) (2.2,2.1) (1.8,2.1) (1.7,2)};
\draw[thick, draw=violet] (1.8,2)--(2.2,2);
\node[violet] at (3,2) {\scriptsize{crack $C$}};
\draw[->,blue] (2,1.6)--(2,1.8);
\node[blue] at (2,1.4) {\scriptsize{transition zone of size $\epsilon$}};
\end{tikzpicture}
\end{minipage}
\begin{minipage}[b]{0.6\textwidth}
\centering
\includegraphics[width=0.6\textwidth]{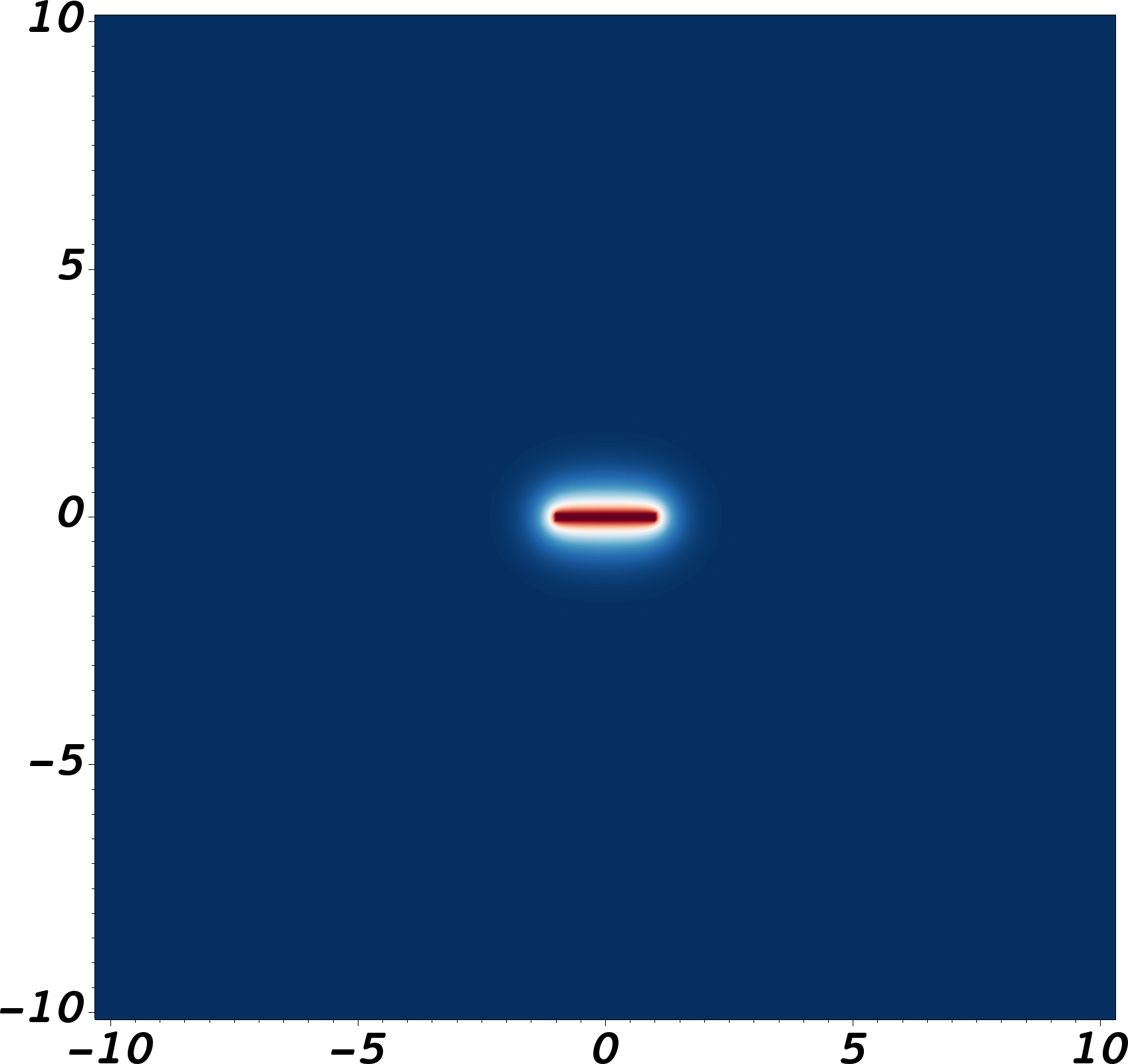}
\end{minipage}
\caption{Domain $\Omega$ (in 2D) with Dirichlet boundaries $\partial\Omega$, an initial crack $C$ of length $2l_0$ and a crack width $\epsilon$, where the phase-field function $\varphi$ is defined.}
\label{domain_Sneddon_2d}
 \end{figure}
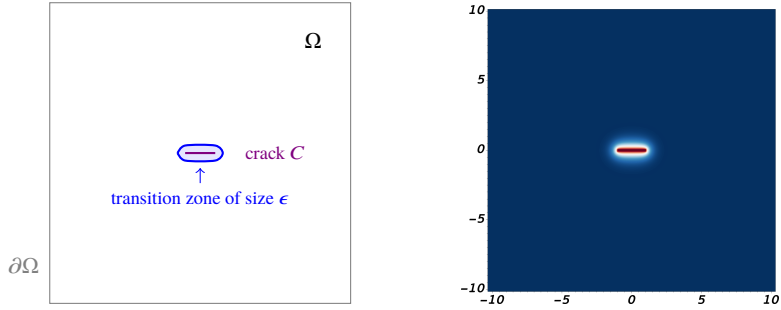
 
Our implementation is based on the open-source software 
DOpElib~\cite{GollWickWollner:2012} and the finite elements from
deal.II~\cite{dealII91,BangerthHartmannKanschat2007}.
The refinement strategy follows~\cite{MangWallothWickWollner:2019}[Section 4.2]. This strategy allows to flag certain cells based on the cell-wise error indicators to reach a grid that is optimal with respect to an objective function that
tries to balance reducing the error and increasing the numerical cost
due to the added unknowns.

%%%%%%%%%%%%%%%%%%%%%%%%%%%%%%%%%%%%%%%%%%%%%%%%%%%%%%
 \paragraph{Setup}
 We follow the setup from~\cite{WiMueKolletal20}, where the case $\nu = 0.2$ is
 discussed.  
 We assume a two-dimensional domain $\Omega = (-10,10)^2$
 as sketched in Figure~\ref{domain_Sneddon_2d}.
 In this domain, an initial crack with length $l_0 = 2.0$ and thickness $d$
 of two cells on $\Omega_c=[-1,1] \times [-d, d] \subset \Omega$ 
 is prescribed by help of the phase-field function $\varphi$, i.e.,
 $\varphi = 0$ in $\Omega_c$ and $\varphi = 1$ in $\Omega\setminus\Omega_c$.
 Note that the thickness of $2d$ corresponds to $2h/\sqrt{2}$, where $h$ is the
 cell diameter.
 For the numerical realization, $\varphi^0_h = I_h^0(\varphi^0)$ is
 utilized.

 As boundary conditions, the displacements $u$ are set to zero on $\partial
 \Omega$.
 For the phase-field variable, we use homogeneous Neumann conditions (so-called traction
 free conditions), i.e., $\epsilon \partial_n \varphi = 0$ on $\partial \Omega$.\\

 For all tests in the following sections the crack bandwidth $\epsilon$ is set
 as $\epsilon = 4\sqrt{2} d$, the regularization parameter $\kappa$ is determined
 sufficiently small with $\kappa =10^{-8}$. 
 The fracture toughness of the observed material is $G_c =1.0$ and the Young's
 modulus $E=1.0$.
 
 The numerical tests in the following are based on three configurations derived from Sneddon's setup as discussed in detail in~\cite{WiMueKolletal20} 
 using the solving strategy described below for the discrete formulations of Section~\ref{sec_disc}
 and adaptively refined meshes based on the error estimator in Section~\ref{sec_error}:
\begin{itemize}
 \item {\bf Example 1}: Constant given pressure with $p_g = 10^{-3}$ and $\nu =0.2$ to $\nu= 0.5$ using
Problem~\ref{form_1_h}, called Example 1A,\\
compared to Problem~\ref{form_1_mixed_h}, called Example 1B, 
in Section~\ref{sneddon};
\item {\bf Example 2}: Constant given pressure with $p_g = 10^{-3}$, $\nu =0.2$ to $\nu= 0.5$ and a compressible layer around the finite domain as well as in the prescribed fracture using
Problem~\ref{form_1_h}, called Example 2A,\\
compared to Problem~\ref{form_1_mixed_h}, called Example 2B, in Section~\ref{sneddon_layered}, where details on the layer will be given;
\item {\bf Example 3}: Non-constant given pressure $p_g$, $\nu= 0.5$ and a compressible layer around the finite domain as well as in the prescribed fracture using
Problem~\ref{form_1_mixed_h}, called Example 3A,\\
compared to Problem~\ref{form_1_mixed_h_rob}, called Example 3B, 
in Section~\ref{sneddon_bump_layered}.
\end{itemize}

\paragraph{Solution algorithm}

The coupled inequality system in
Problems~\eqref{form_1_h},~\ref{form_1_mixed_h}, and~\ref{form_1_mixed_h_rob} is formulated as a complementarity system as shown in~\cite{MaWiWo20}.
Therein a Lagrange multiplier is introduced for treating the inequality constraint.
The Lagrange multiplier $\tau$ is discretized in 
the dual basis to the $\mathbb{Q}_1$ space denoted by $\mathbb{Q}_1^*$
and the corresponding discrete function space denoted as $X_h$.

The discrete form is then solved in a monolithic
fashion, but noticing that 
$\varphi$ is time-lagged in the first term of the displacement
equation. 
This means in Problem~\ref{form_1_h},~\ref{form_1_mixed_h},
and~\ref{form_1_mixed_h_rob}
we replace in~\eqref{eq_u_mixed_disc},~\eqref{eq_u_disc} and~\eqref{eq_u_mixed_disc_rob}, respectively, the term $g(\varphi_h)$ by 
$g(\varphi_h^{n-1})$ and $(\varphi_h^{n})^2$ by $(\varphi_h^{n-1})^2$. This procedure helps in relaxing the nonlinearity. Of course,
a temporal discretization error is introduced, which however is not significant 
in the steady-state tests considered here.
To this end, we 
formulate a compact form by summing up all equations:
Given the initial data $\varphi^0$; 
for the loading steps $n=1,2,\ldots, N$:\\
Find $U_h:= U_h^n = (\bu_h,p_h,\varphi_h,\tau_h) 
\in Y_h:= (V_h\times P_h\times W_h\times X_h)$ such that
\[
A_{\varphi^{n-1}}(\bu_h,p_h,\varphi_h,\tau_h) = 0.
\]
To solve $A_{\varphi^{n-1}}(\cdot) = 0$, we formulate 
a residual-based Newton scheme, e.g.,~\cite{Wi17_SISC}.
The concrete scheme (and its implementation) 
can be found in~\cite{GollWickWollner:2012,dope}.
The occurring linear systems are
solved with a direct method provided by UMFPACK~\cite{DaDu97}.

%%%%%%%%%%%%%%%%%%%%%%%%%%%%%%%%%%%%%%%%%%%%%%%%%%%%%%
\paragraph{Quantities of interest}
For all examples, we compared the following quantities of interest:
\begin{itemize}
\item  Total crack volume (TCV);
\item  Bulk energy $E_b$;
\item  Crack energy $E_c$.
\end{itemize}
It will turn out for the discussion below, that focusing on TCV
will be sufficient.

For the TCV, manufactured reference values 
can be computed for a infinite domain from the formulae presented in
\cite{SneddLow69}[Section 2.4]. Numerical values on the cut-off domain in
Figure~\ref{domain_Sneddon_2d} and $\nu = 0.2$ can be found in~\cite{WiMueKolletal20}.
Numerically, the total crack volume can be computed by using
\begin{align}
\label{eq_TCV_num}
\text{TCV}=\int_{\Omega} u(x,y) \cdot \nabla \varphi(x,y)\ \mathrm{d}{(x,y)}.
\end{align}

Using the exact representation of $u_y$ (cf.~\cite{SneddLow69}, page
29) applied to our parameter settings as in~\cite{WiMueKolletal20}, we consequently
obtain the reference values listed in Table~\ref{tab:TCV-infinite-reference} 
for an infinite domain.

\begin{table}[htbp!]
\centering
\caption{Manufactured reference values of the TCV computed by help of the formula in 
\cite{WiMueKolletal20} for a infinite domain and different Poisson ratios up to the incompressible limit.}
\renewcommand*{\arraystretch}{1.35}
\begin{tabular}{|l|r|} \hline
  $\nu$ & TCV$_{\text{2d}}$ (reference) \\
   \hline
   $0.2$ & $6.03186\times 10^{-3}$ \\
   $0.49$ & $4.77459\times 10^{-3}$\\
   $0.49999$ & $4.71245\times 10^{-3}$\\
   $0.5$ & $4.71239\times 10^{-3}$\\ \hline
\end{tabular}
\label{tab:TCV-infinite-reference}
\end{table}
As a second quantity of interest, the bulk energy $E_b$ is defined as
\begin{align}
 E_b:=\int_{\Omega} \frac{g(\phi)}{2} \sigma : e(\bu) \mathrm{d}x, \label{bulk}
\end{align}
where $\sigma := \sigma(\bu)$ for Problem~\ref{form_1_h} and $\sigma := \sigma(\bu,p)$ 
for Problem~\ref{form_1_mixed_h} and~\eqref{form_1_mixed_h_rob}.
As third quantity of interest, the crack energy is computed via
\begin{align}
  E_c:=\frac{G_c}{2} \int_{\Omega}  \left(\frac{(\varphi-1)^2}{\epsilon} + \epsilon|\nabla \varphi|^2\right) \ \mathrm{d}x.\label{crack}
\end{align}

%%%%%%%%%%%%%%%%%%%%%%%%%%%%%%%%%%%%%%%%%%%%%%%%%%%%%%
\subsection{Sneddon-inspired test cases (Example 1)}\label{sneddon}
In this first set of numerical examples, we compare Example 1A with Example 1B. 
The prescribed pressure is $p_g = 10^{-3}$ and 
the Poisson ratios are $\nu = 0.2, 0.49, 0.49999$ 
and $\nu=0.5$ (only for the mixed formulation Example 1B).

The starting meshes are once globally uniformly refined and three 
times further uniformly refined around the crack. 
The following three meshes are either uniformly refined in a zone
around the crack (geometric refinement) or adaptively refined based on the estimator proposed in Section~\ref{sec_error}.

\begin{table}[htbp!]
\centering
\caption{The number of degrees of freedom (DoF) and the TCV 
  for four different Poisson ratios for Example 1A and Example 1B.}
\renewcommand*{\arraystretch}{1.35}
\begin{tabular}{|l|l||r|r|r|r|r|r|r|r|}\hline
  \multirow{3}{*}{$\nu$} & \multirow{3}{*}{$d$} & \multicolumn{4}{|c|}{Example 1A} & \multicolumn{4}{|c|}{Example 1B}\\ \hline
  & & \multicolumn{2}{|c|}{geometric} &\multicolumn{2}{|c|}{adaptive} &\multicolumn{2}{|c|}{geometric} &\multicolumn{2}{|c|}{adaptive}\\ \hline
   & & DoF & TCV & DoF & TCV & DoF & TCV & DoF & TCV \\
   \hline
   \multirow{4}{*}{$0.2$}
   & 0.0625    & 29,988  & 0.00818  & 29,988& 0.00818 & 96,436    & 0.00821  &96,436 &0.00821 \\
   & 0.03125   & 74,852  & 0.00691  &36,964 & 0.00691 & 241,860   & 0.00693  & 118,916&0.00693 \\ 
   & 0.015625  & 241,156 & 0.00639  &49,044 &0.00639  & 781,604   & 0.00640  &157,900 &0.00640 \\  
   & 0.0078125 & 880,740 & 0.00616  &69,428 &0.00616  & 2,858,788 & 0.00617  &  223,692 &0.00617 \\ \hline
   \multirow{4}{*}{$0.49$}
   & 0.0625    &29,988   & 0.00601  & 29,988& 0.00601& 96,436    & 0.00620  & 96,436&  0.00620 \\ 
   & 0.03125   & 74,852  & 0.00492  &36,580 &0.00491 & 241,860   & 0.00504  & 117,668& 0.00504  \\ 
   & 0.015625  & 241,156 & 0.00440  &48,436 &0.00438 & 781,604   & 0.00448  & 155,948&  0.00447 \\ 
   & 0.0078125 & 880,740 & 0.00415  & 68,628&0.00413 & 2,858,788 & 0.00421  &221,116 & 0.00421  \\ \hline 
   \multirow{4}{*}{$0.49999$}
   & 0.0625    & 29,988   & 2.28E-5 & 29,988  &2.28E-5 & 96,436    & 2.38E-5 &96,436 &  2.38E-5 \\ 
   & 0.03125   & 74,852   & 2.33E-5 &37,332   & 2.29E-5& 241,860   & 2.39E-5 &120,124 & 2.39E-5  \\ 
   & 0.015625  & 241,156  & 2.35E-5 & 47,844  & 2.28E-5& 781,604   & 2.39E-5 & 154,036& 2.39E-5  \\ 
   & 0.0078125 & 880,740  & 2.36E-5  & 70,164 &2.28E-5 & 2,858,788 & 2.39E-5 &226,108 & 2.39E-5  \\\hline
   \multirow{4}{*}{$0.5$}
   & 0.0625    & & & & & 96,436    &  1.44E-5 &96,436  & 1.44E-5 \\ 
   & 0.03125   & & & & & 241,860   &  -1.15E-6&118,292 &-2.03E-7  \\ 
   & 0.015625  & & & & & 781,604   & -1.29E-7 &155,284 & -1.33E-7 \\ 
   & 0.0078125 & & & & & 2,858,788 & -3.38E-8 &227,356 &-3.46E-8  \\\hline
\end{tabular}
\label{tab:TCV-case1}
\end{table}

\begin{table}[htbp!]
\centering
\caption{The number of degrees of freedom (DoF) and the bulk energy $E_b$ 
  for four different Poisson ratios for Example 1A and Example 1B.}
\renewcommand*{\arraystretch}{1.35}
\begin{tabular}{|l|l||r|r|r|r|r|r|r|r|}\hline
  \multirow{3}{*}{$\nu$} & \multirow{3}{*}{$d$} & \multicolumn{4}{|c|}{Example 1A} & \multicolumn{4}{|c|}{Example 1B}\\ \hline
  & & \multicolumn{2}{|c|}{geometric} &\multicolumn{2}{|c|}{adaptive} &\multicolumn{2}{|c|}{geometric} &\multicolumn{2}{|c|}{adaptive}\\ \hline
   & & DoF & $E_b$ & DoF & $E_b$ & DoF & $E_b$ & DoF & $E_b$ \\ \hline
   \hline
   \multirow{4}{*}{$0.2$}
   & 0.0625    & 29,988  & 4.06E-6  &29,988 &4.06E-6 & 96,436     & 4.07E-6 &96,436   & 4.07E-6\\
   & 0.03125   & 74,852  & 3.38E-6  &36,964  &3.38E-6 & 241,860  & 3.39E-6 & 118,916 & 3.39E-6\\ 
   & 0.015625  & 241,156 & 3.14E-6  & 49,044& 3.13E-6& 781,604   & 3.14E-6 &157,900  &3.14E-6 \\  
   & 0.0078125 & 880,740 & 3.04E-6  & 69,428&3.04E-6 & 2,858,788 & 3.05E-6 & 223,692 & 3.05E-6\\ \hline
   \multirow{4}{*}{$0.49$}
   & 0.0625    &29,988   & 3.00E-6  &29,988 &3.00E-6 & 96,436    & 3.09E-6 &96,436   &3.09E-6\\ 
   & 0.03125   & 74,852  & 2.46E-6  & 36,580&2.45E-6 & 241,860   & 2.52E-6  & 117,668&2.52E-6 \\ 
   & 0.015625  & 241,156 & 2.20E-6  &48,436 & 2.19E-6 & 781,604  & 2.23E-6  & 155,948&2.23E-6 \\ 
   & 0.0078125 & 880,740 & 2.08E-6  & 68,628&2.06E-6 & 2,858,788 & 2.10E-6 &221,116  &2.10E-6 \\ \hline 
   \multirow{4}{*}{$0.49999$}
   & 0.0625    & 29,988   & 1.14E-8 &29,988 &1.14E-8 & 96,436     & 1.19E-8  & 96,436 & 1.19E-8\\ 
   & 0.03125   & 74,852  & 1.16E-8 &37,332 &1.14E-8 & 241,860    & 1.19E-8  &120,124 & 1.19E-8\\ 
   & 0.015625  & 241,156 & 1.17E-8 &47,844 &1.14E-8 & 781,604    & 1.19E-8  & 154,036& 1.19E-8\\ 
   & 0.0078125 & 880,740 & 1.18E-8 &70,164 & 1.14E-8 & 2,858,788 & 1.19E-8  & 226,108&1.19E-8 \\\hline
   \multirow{4}{*}{$0.5$}
   & 0.0625    & & & & & 96,436    &  4.06E-9   &96,436  &4.06E-9\\ 
   & 0.03125   & & & & & 241,860   &-3.25E-10   &118,292 &-5.74E-11\\ 
   & 0.015625  & & & & & 781,604   &-3.65E-11   &155,284 &-3.76E-11\\ 
   & 0.0078125 & & & & & 2,858,788 & -9.56E-12  &227,356 &-9.77E-12\\\hline
\end{tabular}
\label{tab:Eb-case1}
\end{table}
 
Tables~\ref{tab:TCV-case1} and~\ref{tab:Eb-case1} show the resulting values for the TCV and $E_b$ on the starting mesh 
and the following three geometrically or adaptively refined meshes with adjusted parameters
$\epsilon$ and $d$ according to~\cite{WiMueKolletal20}. 

\begin{remark}
  Considering adaptively refined meshes, the parameters $\eps$ and $d$ are decreased by a factor of two after
  each refinement.
 Hence these values are the same for the computations on geometrically refined meshes, which allows
 a fair comparison of results coming from geometrically and adaptively refined meshes.
\end{remark}

For $\nu=0.2$ the
TCV and $E_b$ computed with Problem~\ref{form_1_h}, rounded to three significant digits, matches the numbers given
in~\cite{WiMueKolletal20}, hence we conclude the correctness of our
implementation.\\
The fracture energy $E_c$ is identical to the values
in~\cite{WiMueKolletal20}.
On the coarsest mesh this corresponds
to $E_c \approx 2.895$ and on the finest mesh we have $E_c \approx 2.423$. 
As the numbers for $E_c$ are independent of $\nu$ and the chosen formulation,
they are not listed separately.

In the following, we will focus on the behavior of TCV for different Poisson ratios and compare it to the reference values of Table~\ref{tab:TCV-infinite-reference} on an infinite domain. First, we see in Table~\ref{tab:TCV-case1} and Table~\ref{tab:Eb-case1}
that both quantities of interest are numerically stable under mesh refinement.
This shows the robustness of our proposed models and their numerical realization.
Second, we observe that more incompressible materials 
yield smaller values of the TCV much smaller than the predicted values
in Table~\ref{tab:TCV-infinite-reference}. Physically, this is to be expected 
if we think of incompressible material in a closed box, because the 
material cannot move.
Due to the cut-off of the computational domain and the use of an incompressible material, 
no movement can be expected for $\nu \approx 0.5$. 
This led us to suggest the setting of Section~\ref{sneddon_layered} 
where we add an artificial compressible layer around
the (nearly) incompressible domain and inside the prescribed fracture
$(-1,1)\times (-d,d)$.

%%%%%%%%%%%%%%%%%%%%%%%%%%%%%%%%%%%%%%%%%%%%%%%%%%%%%%
\subsection{Incompressible material surrounded with a 
compressible layer (Example 2)}\label{sneddon_layered}
As we have seen in the previous example in terms of the total crack volume, 
for $\nu_s = 0.49999$, the fracture in incompressible 
solids will not open anymore and the TCV is almost $0$. On the other hand, 
the formulae in~\cite{SneddLow69}[Section 2.4] suggest a value greater than zero. 
The reason being that therein an infinite domain was assumed.
To study incompressible solids in larger domains, we use a trick and add a compressible 
layer as surrounding area. Considering Figure~\ref{domain_Sneddon_2d}, now we work in a domain $(-20,20)^2$ which contains the previously defined domain $(-10,10)^2$.
The surrounding layer of width $10$ is defined as a compressible material with $\nu = 0.2$.
All other parameters, namely $E$, $G_c$, $\kappa$ and $\Omega_c$ are
kept as before with the values listed in the first paragraph of Section~\ref{sec_tests}.
The same compressible material is used inside of the prescribed
fracture on the set $(-1,1)\times (-d,d)$.

\begin{figure}[htbp!]
\centering
\hfill\includegraphics[width=0.47\textwidth]{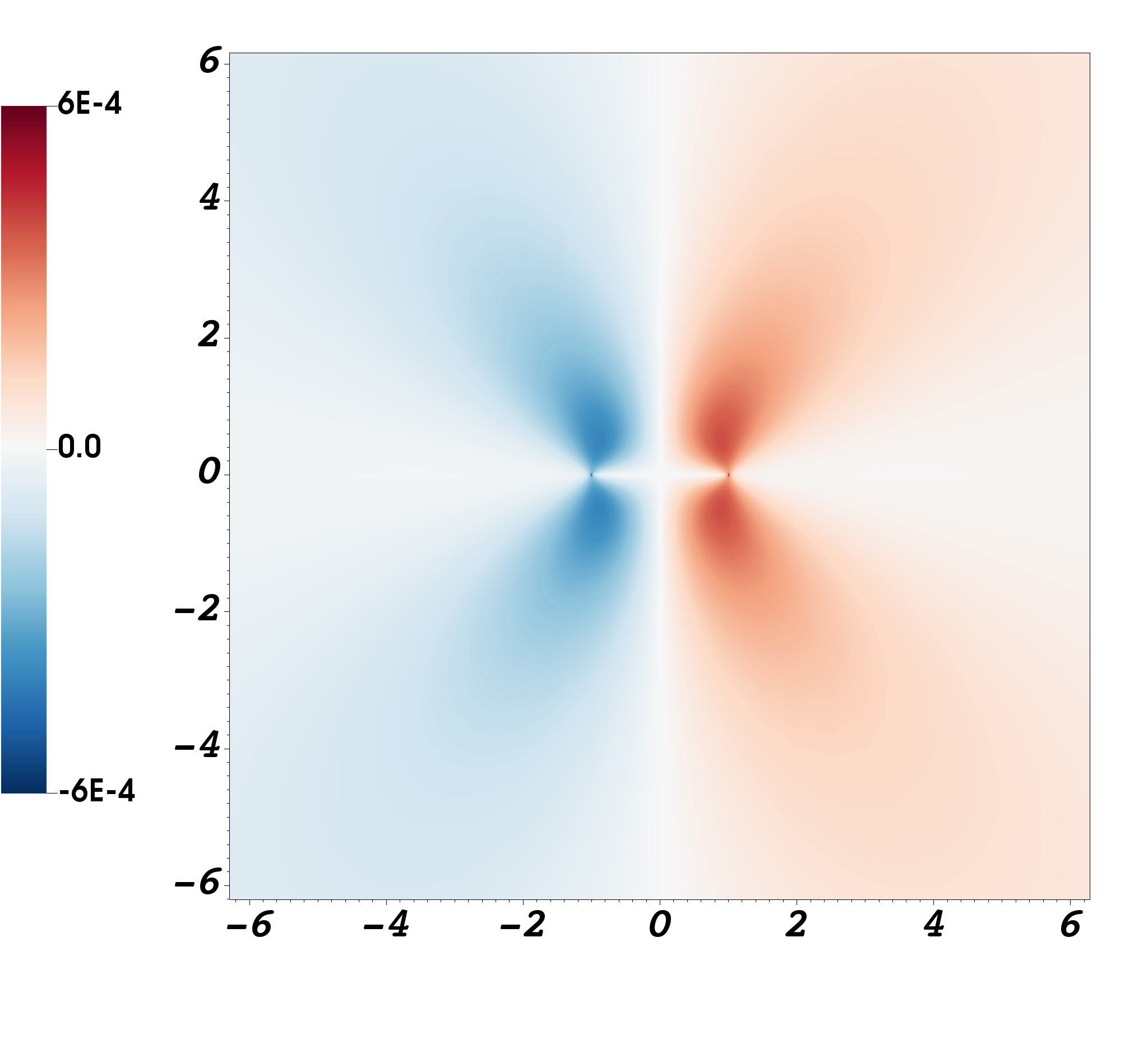}
\hfill\includegraphics[width=0.47\textwidth]{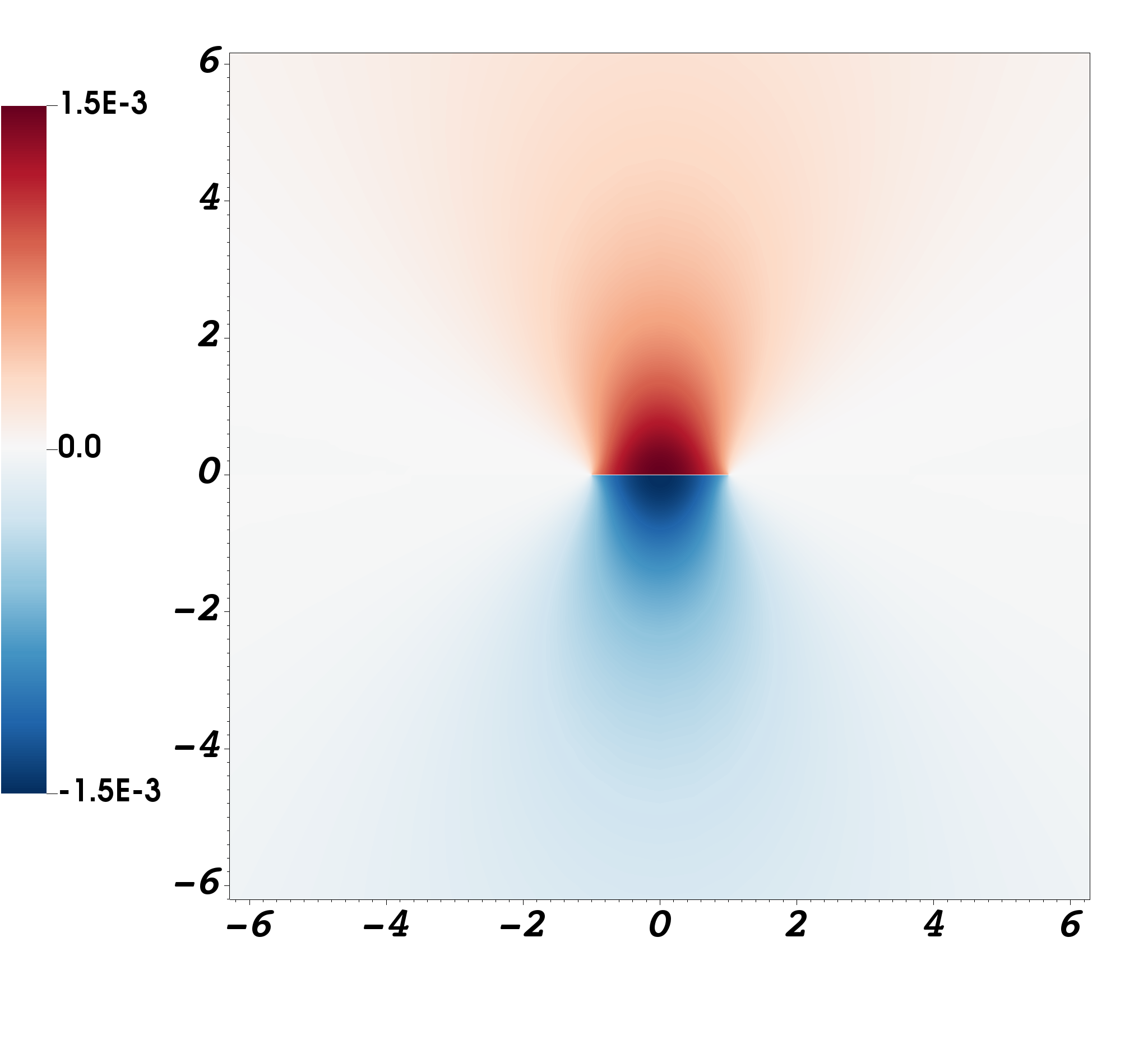}
\hspace*{\fill} \\
\includegraphics[width=0.47\textwidth]{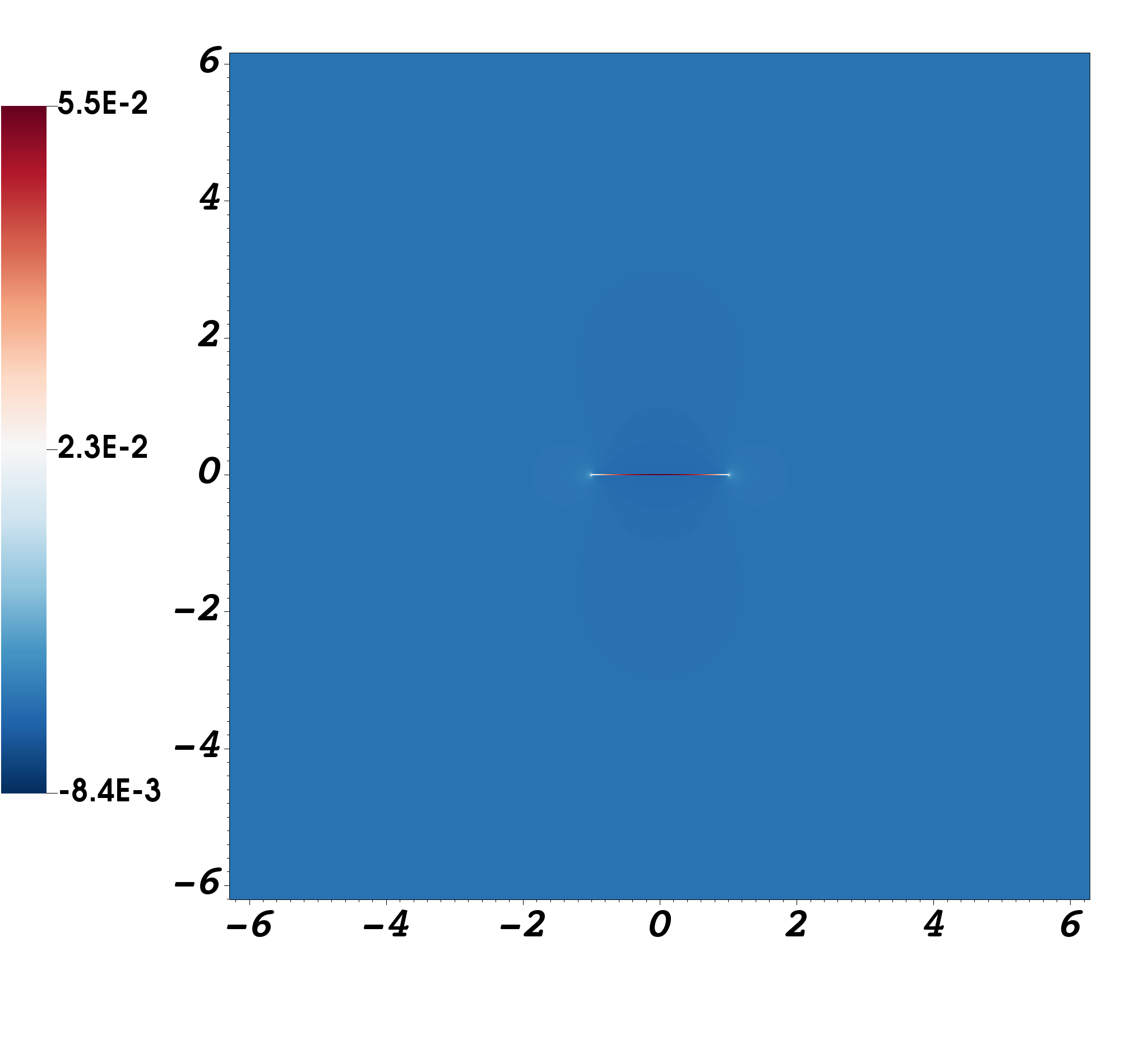}
\caption{Example 2B: The $x$- and $y$-displacements and the pressure $p$ for $\nu=0.5$.
}
\label{fig:x_y_p}
\end{figure}

In Figure~\ref{fig:x_y_p}, the ranges of the $x-$ and the $y-$ displacements as well as for the pressure values are depicted for Example 2B, where a perfect symmetry of the test setup can be observed.

In Table~\ref{tab:TCV-case2}, for the primal-based form (Example 2A), the TCV is underestimated for $\nu \approx 0.5$ while the mixed form (Example 2B) gives results 
consistent with the computations for $\nu = 0.5$. Compared to Table~\ref{tab:TCV-infinite-reference}, the TCV values based on the mixed form (Example 2B) are very similar for the four listed Poisson ratios compared to the 
reference values. Keep in mind at this point, that the reference values are given analytically considering an infinite domain.
\begin{table}[htbp!]
\centering
\caption{The number of degrees of freedom (DoF) and the TCV 
  for four different Poisson ratios for Example 2A and Example 2B.}
\renewcommand*{\arraystretch}{1.35}
\begin{tabular}{|l|l||r|r|r|r|r|r|r|r|}\hline
  \multirow{3}{*}{$\nu$} & \multirow{3}{*}{$d$} & \multicolumn{4}{|c|}{Example 2A} & \multicolumn{4}{|c|}{Example 2B}\\ \hline
  & & \multicolumn{2}{|c|}{geometric} &\multicolumn{2}{|c|}{adaptive} &\multicolumn{2}{|c|}{geometric} &\multicolumn{2}{|c|}{adaptive}\\ \hline
   & & DoF & TCV & DoF & TCV & DoF & TCV & DoF & TCV \\
   \hline
   \multirow{4}{*}{$0.2$}
   & 0.0625    & 49,508  & 0.00836 &49,508 &0.00836 & 159,316   & 0.00839 & 159,316  &0.00839 \\
   & 0.03125   & 94,372  & 0.00703 &58,036 &0.00702 & 304,740   & 0.00704 & 186,828  &0.00704 \\ 
   & 0.015625  & 260,676 & 0.00648 &72,420 &0.00648 & 844,484   & 0.00649 &233,300   &0.00649 \\  
   & 0.0078125 & 900,260 & 0.00624 &93,220 & 0.00624& 2,921,668 & 0.00625  & 300,420  & 0.00625\\ \hline
   \multirow{4}{*}{$0.49$}                                        
   & 0.0625    & 49,508  & 0.00808 &49,508 & 0.00808& 159,316   & 0.00842 & 159,316 &0.00842 \\
   & 0.03125   & 94,372  & 0.00622 &58,420 &0.00620 & 304,740   & 0.00640 &188,076  & 0.00640 \\ 
   & 0.015625  & 260,676 & 0.00540 &72,804 &0.00537 & 844,484   & 0.00551 &234,548  &0.00551 \\  
   & 0.0078125 & 900,260 & 0.00503 &93,796 &0.00500 & 2,921,668 & 0.00511 &301,668  &0.00511 \\ \hline
   \multirow{4}{*}{$0.49999$}                                     
   & 0.0625    & 49,508  & 0.000913 &49,508 & 0.000913 & 159,316   & 0.00840 & 159,316&0.00840 \\
   & 0.03125   & 94,372  & 0.00129 &57,236 &0.000890  & 304,740   & 0.00636 &188,076 &0.00636 \\ 
   & 0.015625  & 260,676 & 0.00188 &71,044 &0.000860  & 844,484   & 0.00545 &234,548 &0.00545 \\  
   & 0.0078125 & 900,260 & 0.00237 &91,220 &0.000838  & 2,921,668 & 0.00505 &301,668 &0.00505 \\ \hline
   \multirow{4}{*}{$0.5$}                                         
   & 0.0625    &         &            & & & 159,316   & 0.00840 &159,316  &0.00840 \\
   & 0.03125   &         &            & & & 304,740   & 0.00636 &188,076  &0.00636 \\ 
   & 0.015625  &         &            & & & 844,484   & 0.00545 & 234,548 &0.00545 \\  
   & 0.0078125 &         &            & & & 2,921,668 & 0.00505 & 301,668 &0.00505 \\ \hline
\end{tabular}
\label{tab:TCV-case2}
\end{table}
 
Further, the TCV in Table~\ref{tab:TCV-case2} on adaptively refined
meshes in comparison to geometrically refined meshes coincide
satisfactorily. Note however, that as it has to be expected the primal
formulation~\eqref{form_1_h} provides unreliable values for $\nu$
close to $0.5$.

\begin{figure}[htbp!]
\centering
\hfill\includegraphics[width=0.45\textwidth]{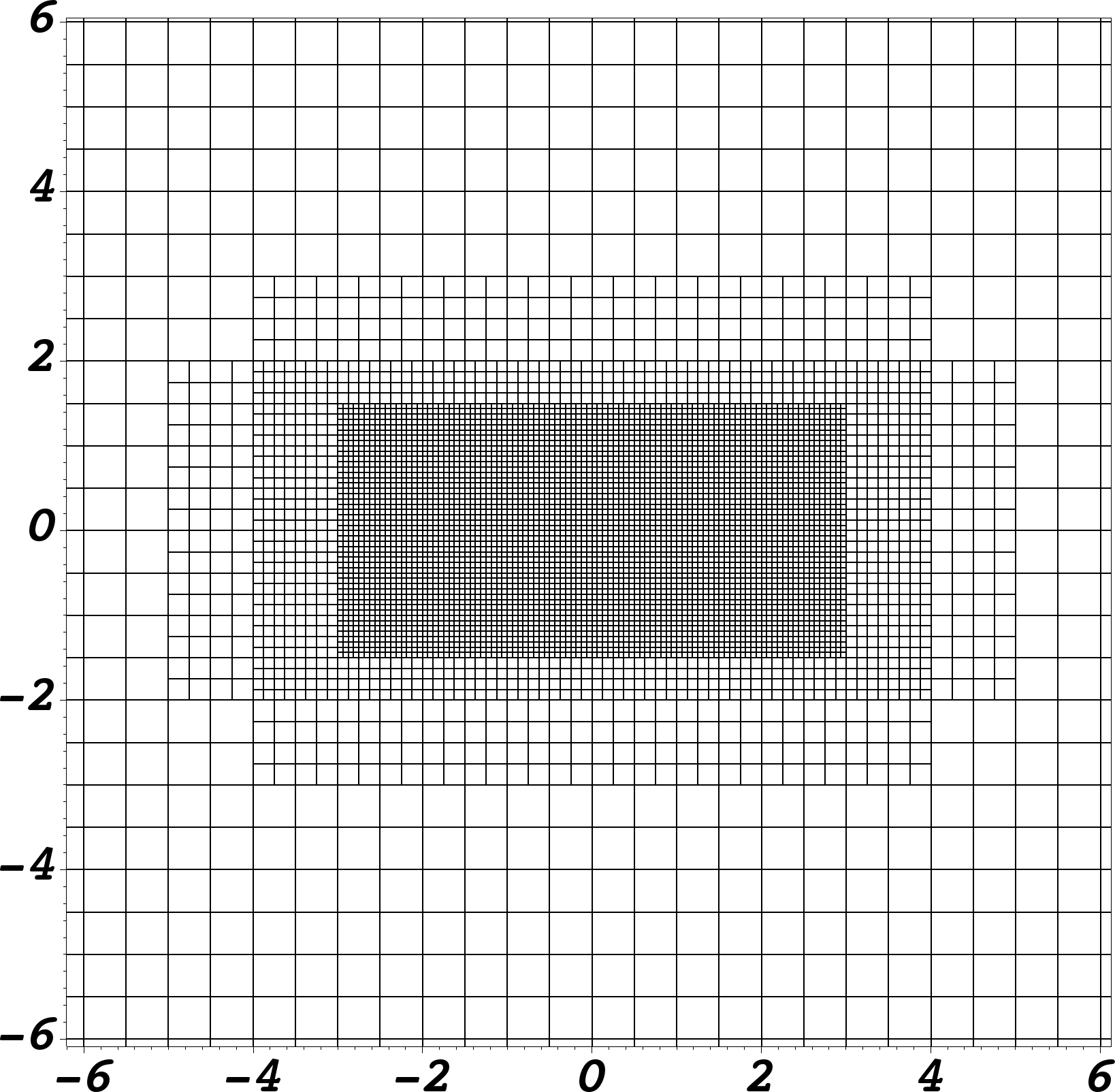}\hfill
\includegraphics[width=0.45\textwidth]{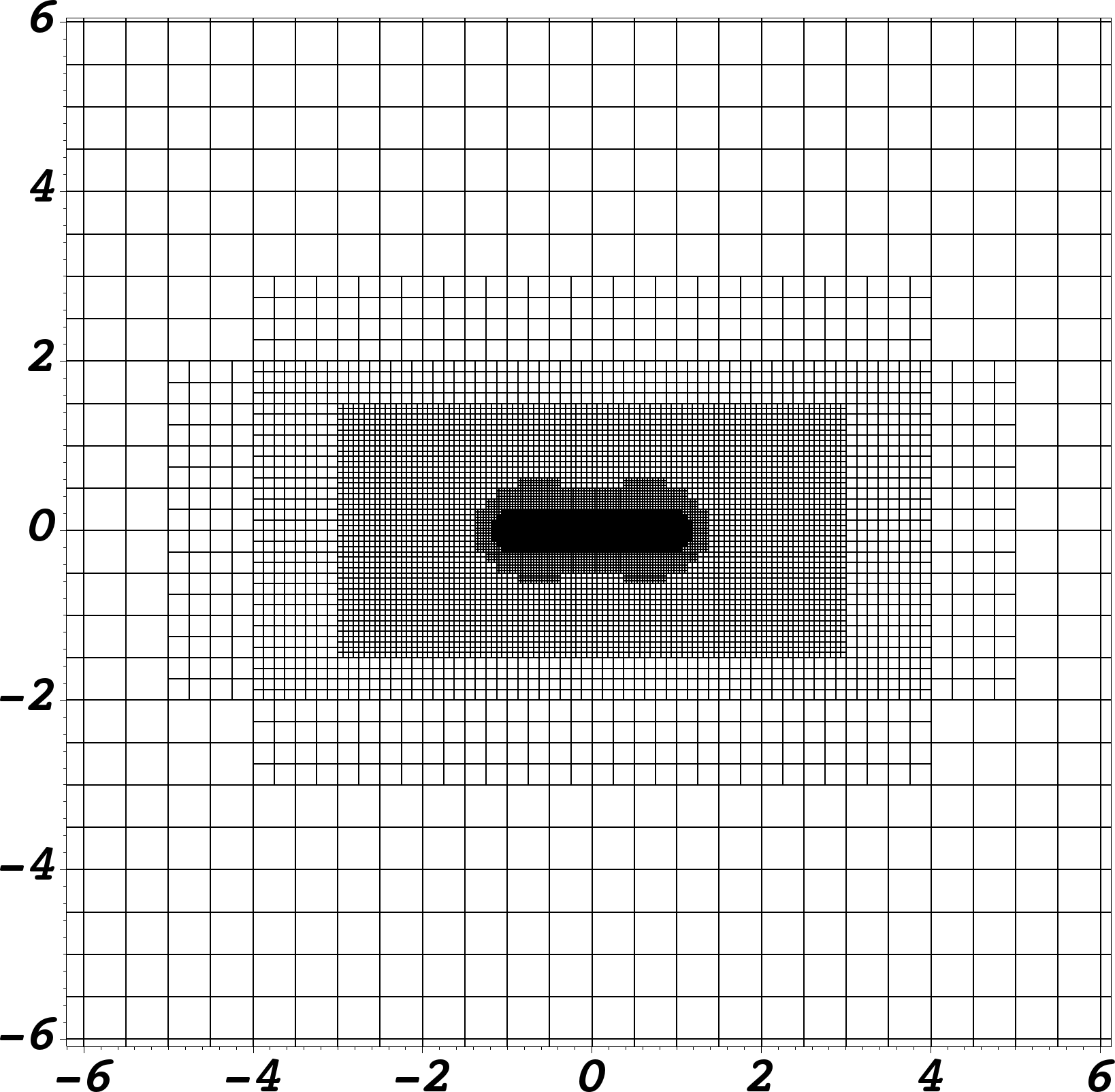}\hspace*{\fill} 
\caption{The mesh on the starting grid (left) and after three levels of adaptive refinement (right) for Example 2B with $\nu=0.5$ zoomed to the crack zone.}
\label{fig:meshes}
\end{figure}

To give an impression of the used meshes and to see the difference between geometrically and adaptively refined meshes, 
in Figure~\ref{fig:meshes}, a coarser starting mesh (geometrically
prerefined) on the left and the mesh after three additional adaptive
refinements  (based on the error  estimator) on the right are given.
Thinking of the problem size and workload, the adaptively refined meshes by help of the error estimator of Section~\ref{sec_error} just needs a tenth of the DoFs, 
but yields very similar results for the TCV on the finest refinement level.

%%%%%%%%%%%%%%%%%%%%%%%%%%%%%%%%%%%%%%%%%%%%%%%%%%%%%%
\subsection{Nonhomogeneous pressure test case with a compressible layer (Example 3)}\label{sneddon_bump_layered}
In this third example, we prescribe a nonhomogeneous pressure $p_g$
in form of a bump that resembles to a fluid-filled fracture situation 
(e.g.,~\cite{MiWheWi15c}). In this situation, we can no longer expect our pressure $p$ to be
almost constant. As it has been observed, e.g., in~\cite{LinkeMatthiesTobiska:2016,LinkeMerdonWollner:2015}
for Stokes flow, for incompressible situations the difficulty in approximating the pressure can
negatively influence the approximation of the displacement field. Hence, for the third example, we will
focus on the case $\nu = 0.5$ and compare the numerical results from Problem~\ref{form_1_mixed_h} with
the pressure robust Problem~\ref{form_1_mixed_h_rob}.

For this setting we consider the following given pressure:
%bump3
\[
p_g(x,y) = f(x)g(y)
\]
where
\begin{align*}
f(x) &= \begin{cases}
  0.001 & 1 \le x < 2, \\
  -0.002\, x^2(x-1.5) & 0 \le x < 1,\\
  0.002\, (x-3)^2(x-1.5) & 2 \le x < 3,\\
  0 & \text{otherwise,}\\
\end{cases}
\\
g(y) &= \begin{cases}
  1 & |y| < 0.5, \\
  2(|y|-1.5)^2|y| & 0.5 \le |y| < 1.5,\\
  0 & \text{otherwise.}\\
\end{cases}
\end{align*}
All other parameters are chosen as in Example~2.

\begin{figure}[htbp!]
\centering
\includegraphics[width=0.49\textwidth]{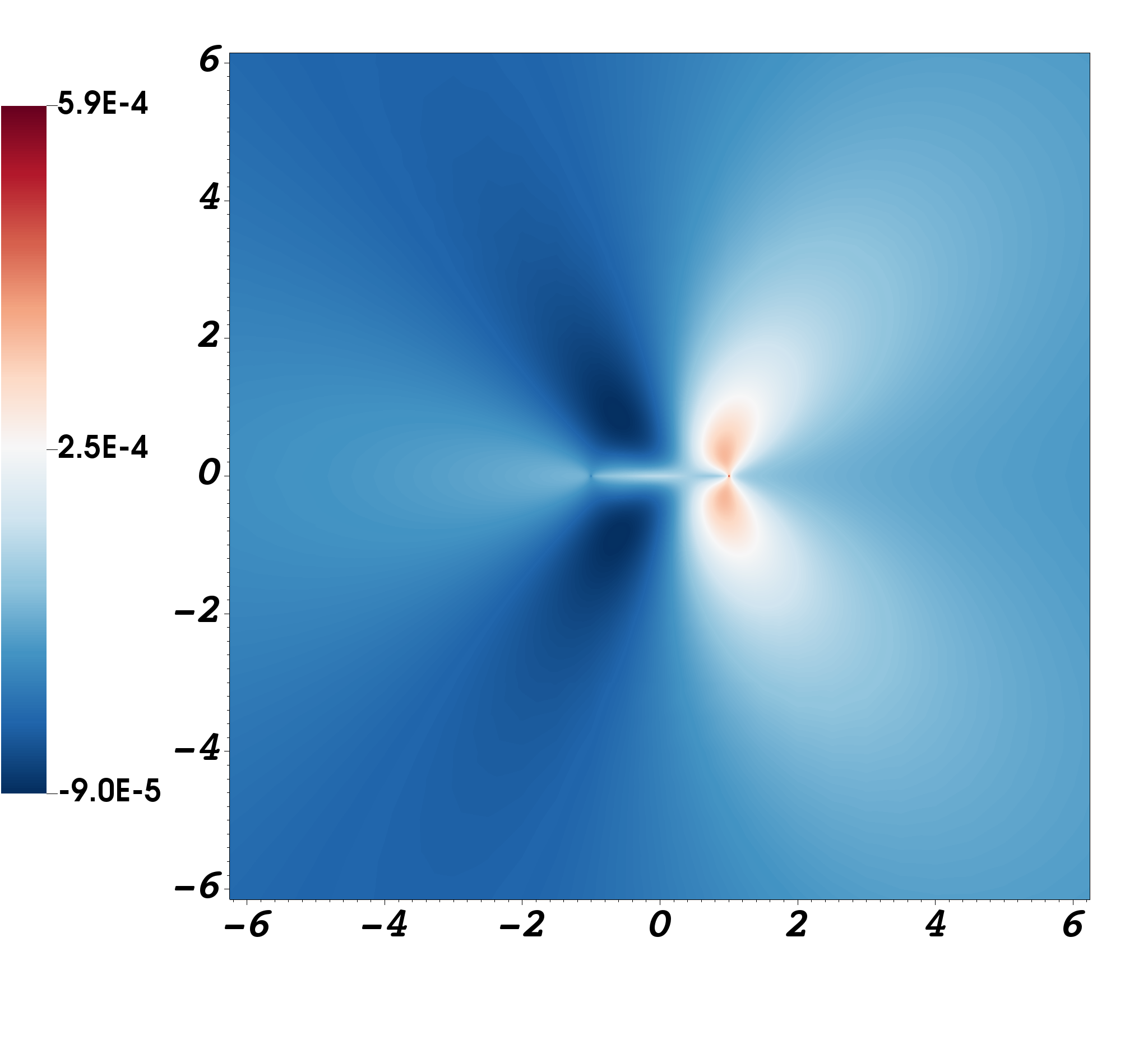}
\includegraphics[width=0.49\textwidth]{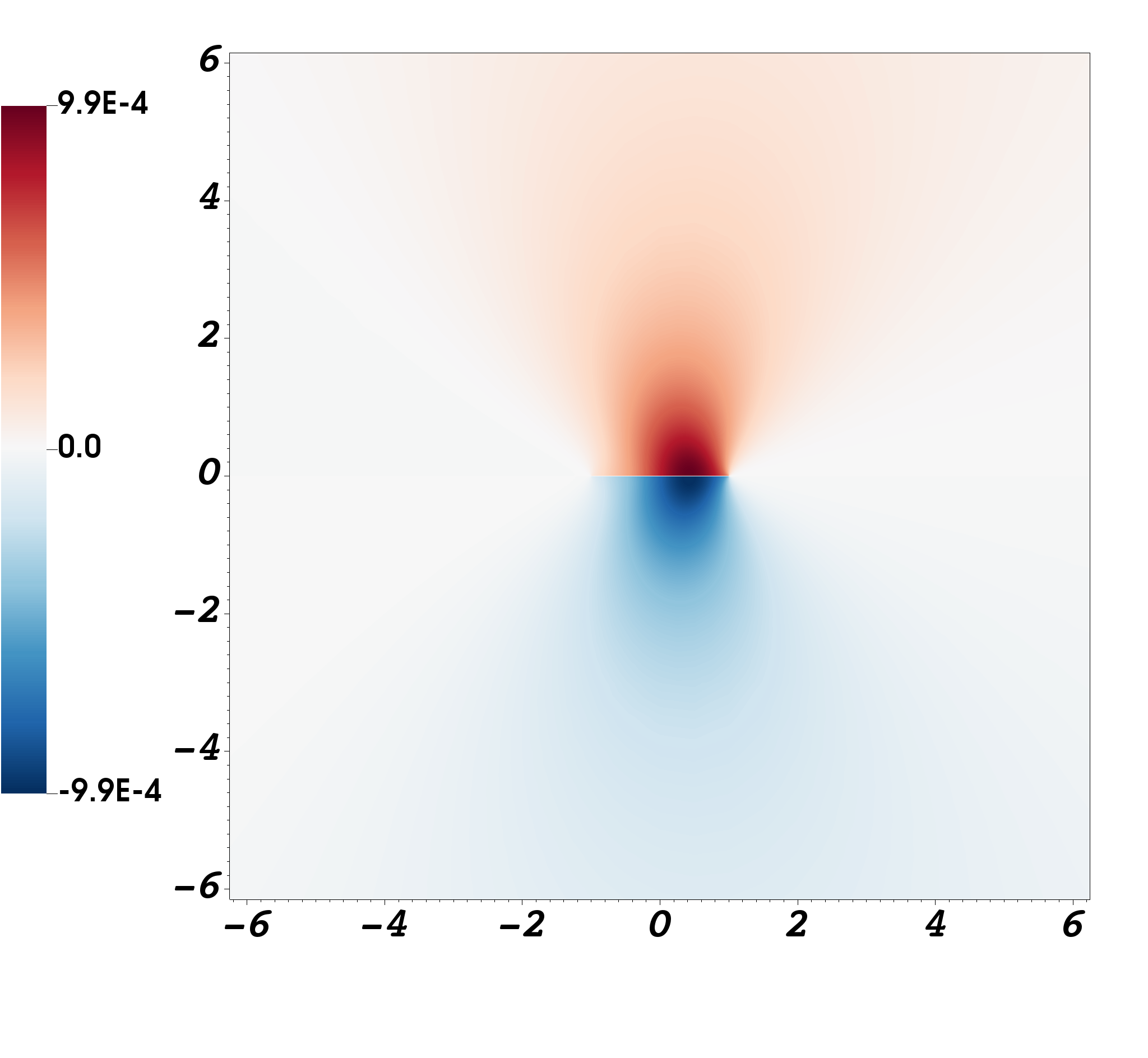}
\includegraphics[width=0.49\textwidth]{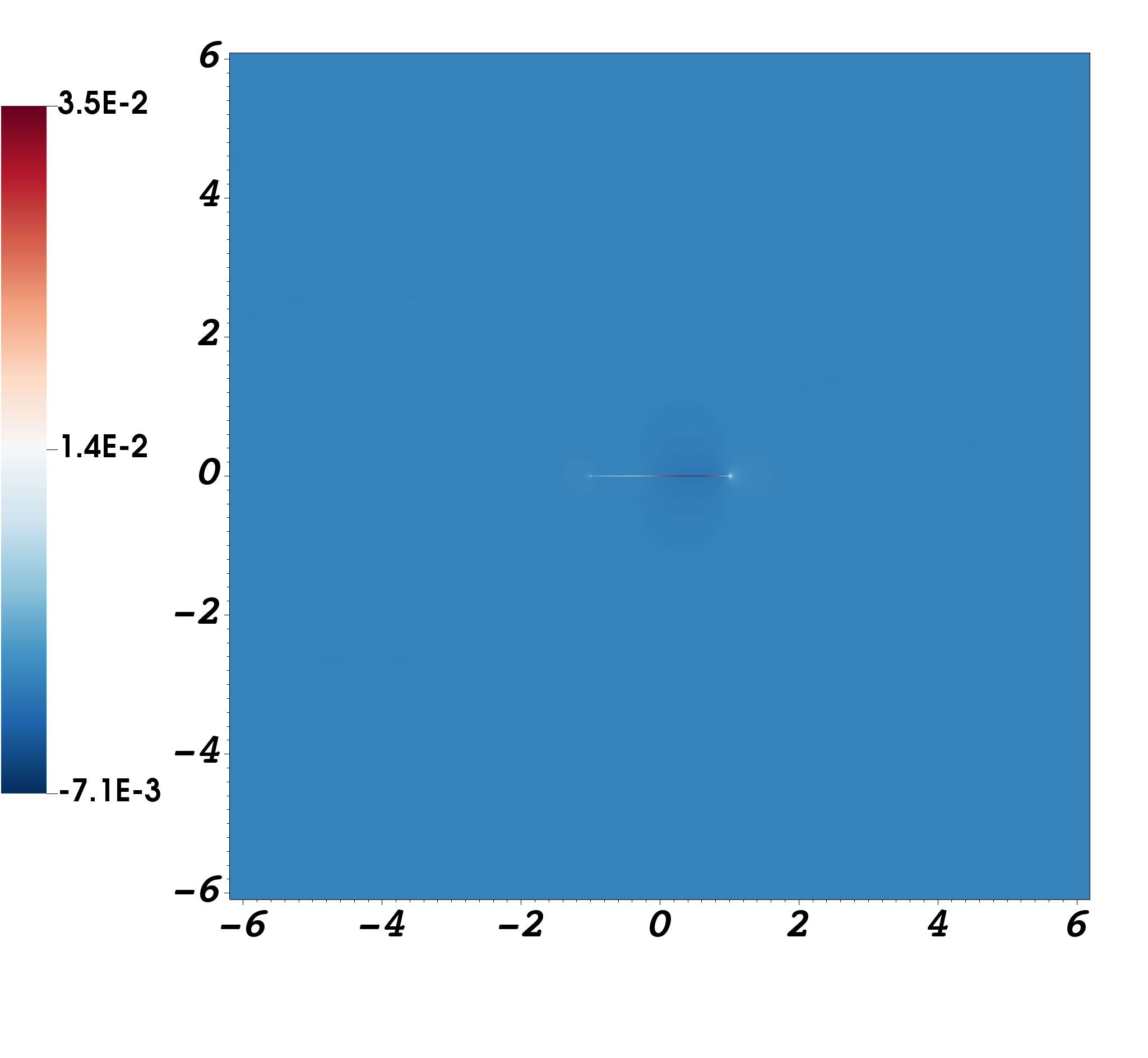}
\includegraphics[width=0.49\textwidth]{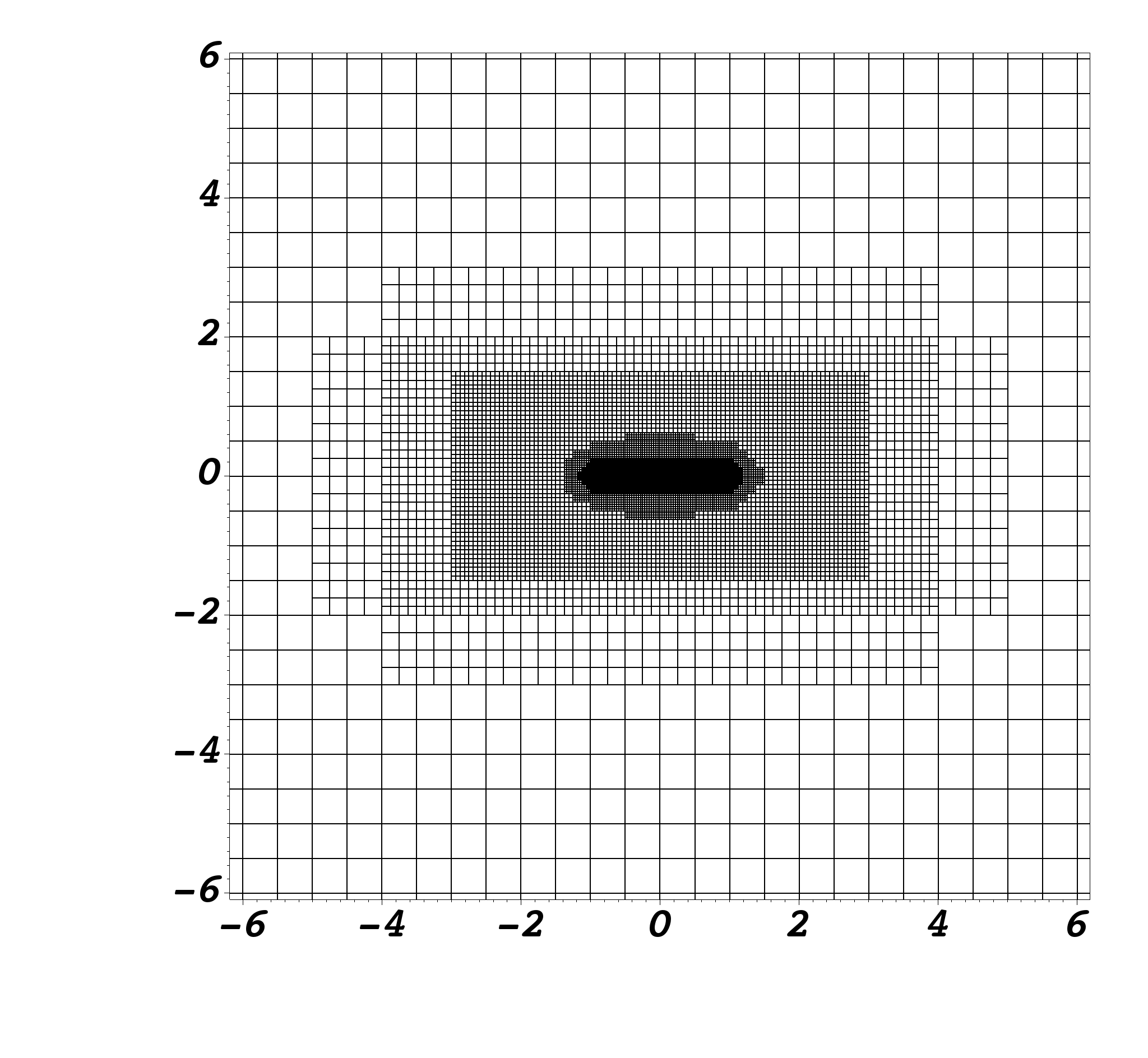}
\caption{Example 3: The $x$- and $y$-displacements (top row), 
and the pressure $p$ for $\nu=0.5$ and the final locally adapted mesh.
}
\label{fig:ex_3}
\end{figure}

The solution is shown in Figure~\ref{fig:ex_3}, where the nonsymmetry
in the setup can be clearly seen in the $x$-displacements. 
It should be noted that similar to Example~2 the pressure is
relatively simple, and the jump in the pressure on the prescribed
fracture is aligned with the mesh. Hence, no difficulty in the
pressure approximation is expected - and thus the pressure robust
results should not deviate too much. Indeed, as the numbers in Table~\ref{tab:TCV-case3}
show the pressure robust discretization yields similar numerical
results. This is not visible in the table, but actual numbers differ in later digits. 
\begin{table}[htbp!]
\centering
\caption{The number of degrees of freedom (DoF) and the TCV 
  for Example 3A and Example 3B.}
\renewcommand*{\arraystretch}{1.35}
\begin{tabular}{|l|l||r|r|r|r|r|r|r|r|}\hline
  \multirow{3}{*}{$\nu$} & \multirow{3}{*}{$d$} & \multicolumn{4}{|c|}{Example 3A} & \multicolumn{4}{|c|}{Example 3B}\\ \hline
  & & \multicolumn{2}{|c|}{geometric} &\multicolumn{2}{|c|}{adaptive} &\multicolumn{2}{|c|}{geometric} &\multicolumn{2}{|c|}{adaptive}\\ \hline
   & & DoF & TCV & DoF & TCV & DoF & TCV & DoF & TCV \\
   \hline
   \multirow{4}{*}{$0.5$}
   & 0.0625    & 159,316   &0.00372 &159,316 &0.00372  & 159,316   & 0.00372 &159,316 & 0.00372 \\
   & 0.03125   & 304,740   &0.00314 &187,744 &0.00314  & 304,740   & 0.00314 &187,744 & 0.00314 \\ 
   & 0.015625  & 844,484   &0.00273 &233,280 &0.00273  & 844,484   & 0.00273 &233,280 & 0.00273 \\  
   & 0.0078125 & 2,921,668 &0.00252 &299,092 &0.00252  & 2,921,668 & 0.00252 &299,092 & 0.00252 \\ \hline
\end{tabular}
\label{tab:TCV-case3}
\end{table}

Since the given pressure only enters the equation on the boundary of
the approximate fracture, i.e., the region where $\nabla \phi \ne 0$,
this rather similar behavior of Problem~\ref{form_1_mixed_h}
and~\ref{form_1_mixed_h_rob} has to be expected. It remains subject to
future research if this remains the same for growing fractures of other
forcings.

%%%%%%%%%%%%%%%%%%%%%%%%%%%%%%%%%%%%%%%%%%%%%%%%%%%%%%
\section{Conclusions}
\label{sec_conclusions}
In this work, we developed a pressurized phase-field fracture model in mixed form for solids up to the incompressible limit $\nu=0.5$. In addition, a residual-type error estimator is presented
for the variational inequality, in this context especially for fractures in solids which are (nearly) incompressible. Estimating the error in the phase-field variable allows to obtain a good resolution especially of the
fracture zone. \\
We investigated the performance of the mixed phase-field fracture formulation and the error estimator with the help of three numerical configurations, all based on Sneddon's and Lowengrub's setup~\cite{SneddLow69} and~\cite{sneddon1946distribution}. 
The theoretical calculations therein based on an infinite pressure-driven cavity and in particular an exact formula for the total crack volume, 
in this work is mainly used to prove the quality of the mixed form as well as the adaptive refinement based on the error estimator.\\
In a second numerical configuration we added a compressible layer around the (nearly) incompressible cavity to allow computing similar results for the TCV as given by the exact formula on an infinite domain.
The findings observed on a compressible layered cavity, which is incompressible in the inner square and around the crack zone, are very convincing. 
To go even further, as a third numerical example, we added a non-constant pressure to the layered Sneddon configuration to provide results of a configuration which is not totally symmetric and tested the results in comparison with a pressure robust modification. It turned out that in the benchmark setup the pressure approximation has no
significant influence on the displacement fields and thus a pressure robust discretization is not necessary.
It will be subject to further studies to check if the situation
remains similar considering a fracture which is not only opening in width but also growing in length.

\begin{acknowledgement}
  Funded by the Deutsche Forschungsgemeinschaft (DFG, German Research
Foundation) -- Projektnummer 392587580
\end{acknowledgement}
%%%%%%%%%%%%%%%%%%%%%%%%%%%%%%%%%%%%%%%%%%%%%%%%%%%%%% 
\bibliographystyle{siam}
\bibliography{./lit}
\end{document}